\documentclass[twoside,a4paper]{amsart}

\usepackage{enumerate}
\usepackage[all]{xy}
\SelectTips{cm}{}
\SilentMatrices

\numberwithin{equation}{section}
\theoremstyle{plain}
\newtheorem{lemma}[subsection]{Lemma}
\newtheorem{theorem}[subsection]{Theorem}

\newtheorem{corollary}[subsection]{Corollary}
\newtheorem{proposition}[subsection]{Proposition}
\theoremstyle{definition}
\newtheorem{construction}[subsection]{Construction}
\newtheorem{definition}[subsection]{Definition}
\newtheorem{example}[subsection]{Example}
\newtheorem{remark}[subsection]{Remark}

\newcommand{\mF}{{\mathbb F}}

\newcommand{\mZ}{{\mathbb Z}}

\newcommand{\cC}{{\mathcal C}}
\newcommand{\cD}{{\mathcal D}}
\newcommand{\cF}{{\mathcal F}}
\newcommand{\cH}{{\mathcal H}}

\newcommand{\cT}{{\mathcal T}}
\newcommand{\cU}{{\mathcal U}}

\newcommand{\cW}{{\mathcal W}}

\DeclareMathOperator{\id}{id}
\DeclareMathOperator{\Hom}{Hom}
\DeclareMathOperator{\Aut}{Aut}

\DeclareMathOperator{\HML}{HML}
\DeclareMathOperator{\Ext}{Ext}
\DeclareMathOperator{\HH}{HH}
\DeclareMathOperator{\THH}{THH}
\DeclareMathOperator{\skel}{sk}
\DeclareMathOperator{\GL}{GL}
\DeclareMathOperator{\Mat}{Mat}
\DeclareMathOperator{\Sq}{Sq}

\DeclareMathOperator{\Objec}{Ob}
\DeclareMathOperator{\Map}{Map}
\DeclareMathOperator*{\colim}{colim}
\DeclareMathOperator{\Ho}{Ho}

\newcommand{\tensor}{\otimes}
\newcommand{\ovl}{\overline}
\newcommand{\udl}{\underline}
\newcommand{\ot}{\leftarrow}
\newcommand{\sm}{\wedge}
\newcommand{\wdg}{\vee}
\newcommand{\iso}{\cong}

\newcommand{\op}{{\mathrm{op}}}

\newcommand{\Ab}{{\mathcal Ab}}
\newcommand{\Top}{{\mathcal Top}}
\newcommand{\cAut}{{\mathcal Aut}}
\newcommand{\Sph}{{\mathbb S}}
\newcommand{\tbracket}[1]{\langle #1 \rangle}
\newcommand{\spl}{{\mathrm{-sp}}}
\newcommand{\modu}[1]{\textrm{Mod-}#1}
\newcommand{\udlhat}[2]{\widehat{\udl{#1}}\textrm{$^{#2}$}}
\newcommand{\ucat}[2]{({#1}\hspace{-.2ex}\downarrow\hspace{-.2ex}{#2})}
\newcommand{\cub}[1]{I^{#1}}
\newcommand{\randcub}[1]{\widetilde{\partial \cub{#1}}}
\newcommand{\tilcub}[1]{\widetilde{\cub{#1}}}
\newcommand{\MapX}{\Map_{\cC}(X,X)}
\newcommand{\MapXinv}{\Map_{\cC}(X,X)^{\times}}
\newcommand{\randsimpl}[1]{\partial \Delta^{#1}}
\newcommand{\BMapX}[1]{B\Map_{\cC}(X^{#1},X^{#1})^{\times}}

\begin{document}
\title{Universal Toda brackets of ring spectra}
\author{Steffen Sagave}
\date{\today}
\address{Department of Mathematics, University of Oslo, Box 1053,  N-0316 Oslo, Norway}
\email{sagave@math.uio.no}
\subjclass[2000]{Primary 55P43; Secondary 19D55, 55S35, 55U35}
\begin{abstract}
We construct and examine the {\em universal Toda bracket} of a highly structured ring spectrum $R$. 
This invariant of $R$ is a 
cohomology class in the Mac Lane cohomology of the graded ring of homotopy groups of $R$
which carries information about $R$ and the category of $R$-module spectra. 
It determines for example all triple Toda brackets of $R$ 
and the first obstruction to realizing a module over the homotopy groups 
of $R$ by an $R$-module spectrum. 

For periodic ring spectra, we study the corresponding theory of {\em higher} universal Toda brackets.
The real and complex $K$-theory spectra serve as our main examples. 
\end{abstract}
\maketitle

\section{Introduction}
In this paper, we study a question about highly structured ring spectra. More
specifically, we construct a cohomological invariant $\gamma_R$ of a ring spectrum $R$, called its
universal Toda bracket, and examine which information about $R$ is encoded in \nolinebreak $\gamma_R$.

We use the term {\em ring spectrum} for what is called an $S$-algebra in \cite{Elmendorf_K_M_M_1997},
a symmetric ring spectrum in \cite{Hovey_S_S_2000}, or an orthogonal ring spectrum in 
\cite{Mandell_M_S_S_2001}. 
A ring spectrum $R$ has an associated module category $\modu{R}$, which 
is a stable model category and has a triangulated homotopy category $\Ho(\modu{R})$. 

For an object $X$ of $\Ho(\modu{R})$, its stable homotopy groups $\pi_*(X)$
form a graded $\pi_*(R)$-module.  
One of our aims is to understand the resulting functor $\pi_*(-) \colon \Ho(\modu{R}) \to \modu{\pi_*(R)}$ better. 
Particularly, we want 
to examine under which conditions a $\pi_*(R)$-module $M$ is {\em realizable}, 
that is, arises as the homotopy groups of an $R$-module spectrum. 

There is an obstruction theory associated to this problem, with obstructions
$\kappa_i(M) \in \Ext^{i,2-i}_{\pi_*(R)}(M,M)$ for $i \geq 3$. The first obstruction $\kappa_3(M)$ is always
defined and unique. It vanishes if and only if $M$ is a retract of a realizable module. For $i \geq 4$,
$\kappa_{i}(M)$ is only defined if $\kappa_{i-1}(M)$ vanishes, and there are choices involved. 
We examine these obstructions and show how they depend on the structure of $R$.

The obstruction theory is the special case of an obstruction theory for realizability in a triangulated category 
$\cT$ described in \cite[Appendix A]{Benson_K_S_04}. In this generality, it 
can be used to find out whether a module over the graded  endomorphism ring $\cT(N,N)_*$ of 
a compact object $N$ can be realized as $\cT(N,X)_*$ for some object $X$ of $\cT$.  
An algebraic instance
of this problem is to realize a module over the cohomology of a differential graded algebra $A$ as the 
cohomology of a differential graded $A$-module.

Because of this analogy between ring spectra and differential graded algebras, the following
result is a motivation for our work: for a differential graded algebra $A$ over a field $k$,
Benson, Krause, and Schwede \cite{Benson_K_S_04} study a class $\gamma_A \in \HH^{3,-1}_{k}(H^*(A))$ in the Hochschild
cohomology of the cohomology ring of $A$. It determines by evaluation all triple (matric) Massey
products of $H^*(A)$. Moreover, via the map 
\[ - \tensor^L \id_M \colon \HH^{3,-1}_{k}(H^* (A)) \to \Ext^{3,-1}_{H^*(A)} (M,M),\]
it determines the first realizability obstruction $\kappa_3(M)$ for every $H^*(A)$-module $M$.

We develop a similar theory for ring spectra. Though the obstruction theory for the realizability problem
takes place completely in triangulated categories, the definition of a cohomology class with that property needs
information from an underlying `model'. In the case of the differential graded algebra $A$, the $A_{\infty}$-structure
of $H^*(A)$ can be used to define $\gamma_A$ \cite[Remark 7.4]{Benson_K_S_04}. In the case of ring spectra, there is no such $A_{\infty}$-structure.
The appropriate replacement is to use 
that choosing representatives in the model category of maps in the homotopy category is in general not associative
with respect to the composition. This non-associativity leads to obstructions which assemble to a well defined cohomology class.

The formulation of our main results uses Mac Lane cohomology groups, denoted by $\HML$. We define this cohomology theory
for graded rings using the normalized cohomology of categories \cite{Baues_W_85}. Its ungraded version is equivalent to 
Mac Lane's original definition \cite{Jibladze_P_1991}. This theory is, for various reasons,
an appropriate replacement of the Hochschild cohomology used in \cite{Benson_K_S_04}. One reason
is that one can, similar to Hochschild cohomology, evaluate a representing cocycle on a sequence of 
composable maps. If the sequence of maps is a complex, it makes sense to ask the evaluation to be an element
of the Toda bracket of the complex. 

One main result is the following special case of  Theorem \ref{univ_higher_tbracket}:
\begin{theorem} \label{triple_universal}
Let $R$ be a ring spectrum. Then there exists a well defined cohomology class  
$\gamma_R \in \HML^{3,-1}(\pi_*(R))$ which, by evaluation, determines
all triple matric Toda brackets of $\pi_*(R)$. For a 
$\pi_*(R)$-module $M$ which admits a resolution by finitely generated free $\pi_*(R)$-modules, 
the product $\id_M \cup \gamma_R \in \Ext^{3,-1}_{\pi_*(R)}(M,M)$
is the first realizability obstruction $\kappa_3(M)$.
\end{theorem}
The term {\em universal Toda bracket} for such a cohomology class, 
as well as the usage of cohomology of categories, are motivated 
by Baues' study of universal Toda brackets for subcategories of the homotopy category of topological spaces 
\cite{Baues_97,Baues_D_89}. The recent preprint \cite{Baues_M_2006} is concerned with a class similar to the
$\gamma_R$ of the last theorem, but studies different properties, namely a relation to ``quadratic pair algebras''. 

Theorem \ref{triple_universal} applies for example to the real $K$-theory spectrum $KO$. As $KO$ has non-vanishing
triple Toda brackets, $\gamma_{KO}$ is non-trivial. Moreover, the obstructions determined
by $\gamma_{KO}$ detect the non-realizable $\pi_*(KO)$-module $(\pi_*(KO)) \tensor \mZ/2$. 
We discuss in Remark \ref{wolberts_thm} how  
this contradicts a claim of Wolbert \cite[Theorems 20 and 21]{Wolbert_1998}.

The proof of Theorem \ref{triple_universal} divides into two parts. In Section \ref{constr_section}, 
we give a general construction of the universal
Toda bracket of a small subcategory of the homotopy category of a stable model category. Specializing to the 
subcategory of finitely generated free modules in $\Ho(\modu{R})$, this defines $\gamma_R$. Theorem \ref{toda_b-obstr} 
shows how a cohomology class which determines Toda brackets also determines the obstructions \nolinebreak $\kappa_3$.

Many examples of ring spectra have the property that their ring of homotopy groups is concentrated in degrees divisible
by $n$ for some $n \geq 2$. Then all realizability obstructions $\kappa_3$ vanish for degree reasons. The first
realizability obstruction not vanishing for degree reasons is determined by a 
{\em higher} universal Toda bracket, which we also introduce in Theorem \ref{univ_higher_tbracket}.

The higher universal Toda bracket of a ring spectrum $R$ becomes particularly nice if $\pi_*(R)$ is
a graded Laurent polynomial ring on a central generator of degree $n$. 
In this case, the higher universal Toda bracket $\gamma_R^{n+2}$ can be defined as an element of $\HML^{n+2}(\pi_0(R))$.
As in Theorem \ref{triple_universal}, 
it determines $(n+2)$-fold Toda brackets and realizability obstructions $\kappa_{n+2}$. 
But now there is a better chance to actually identify $\gamma_R^{n+2}$, since 
computations of (ungraded) Mac Lane-cohomology groups are known in relevant cases.
For example, the universal Toda bracket $\gamma_{KU}^4$ of the complex $K$-theory spectrum $KU$ is an element
of $\HML^4(\mZ) \iso \mZ/2$. We prove in Proposition \ref{computing_gamma_ku} that it is the non-zero element. 

The calculation of $\gamma^4_{KU}$ is a consequence of a different kind of information 
detected by universal Toda brackets. The Toda brackets of a ring spectrum $R$ can be considered as higher order information 
about zero divisors in $\pi_*(R)$ and its matrix rings. It turns out that the universal Toda bracket 
also knows about the units of $R$ and its matrix rings. 

To make the slogan precise, recall that for a ring spectrum 
$R$ and $q \geq 1$, there is a path connected space $B \GL_q R$. It is the classifying space of the topological monoid
given by the invertible path components of the mapping space $\Map_{R}(R^q,R^q)$. 
The algebraic $K$-theory of $R$ can be built from the spaces $B \GL_q R$ \cite{Waldhausen_1978,Elmendorf_K_M_M_1997}.

If $\pi_*(R)$ is concentrated in degrees divisible by $n$ for some $n\geq 1$, we know that $\pi_k (B \GL_q R) = 0$ for
$1 < k < n+1$.
The following corollary follows from Corollary \ref{univ_tbracket_gr_vs_ugr} and Theorem \ref{gammaR_kinvariants}.
It is related  to \cite{Igusa_1982}, see also \cite[Example 4.9, Theorem 3.10]{Baues_D_89}. 

\begin{corollary}
Let $R$ be a ring spectrum such that $\pi_*(R)$ is a Laurent polynomial ring on a central generator in degree $n$. The
restriction map 
\[ \HML^{n+2}(\pi_0(R)) \to H^{n+2}(\pi_1(B \GL_q R), \pi_{n+1}(B \GL_q R)) \]
sends $\gamma_{R}^{n+2}$ to the first $k$-invariant of $B \GL_q R$ not vanishing for dimensional reasons. 
\end{corollary}
Moreover, with an additional assumption on $\HML^{n+1}(\pi_0(R))$, we interpret the vanishing
of $\gamma_R^{n+2}$ in terms of  algebraic $K$-theory in Proposition \ref{splitting_k-theory}.

\subsection*{Organization}
The main results can be found in Section 8. 
There we also discuss the examples mentioned in the introduction.

In the second section, we briefly review cohomology of categories
and Mac Lane cohomology, including a version for graded rings, and define the cup product used in Theorem \ref{triple_universal}. 
In the third section, we explain the obstruction theory for 
realizability in triangulated categories. The fourth section is devoted to (higher) Toda brackets
in triangulated categories. We explain how Toda brackets determine realizability obstructions.   

The fifth section is the technical backbone of this paper. We give a general construction of the universal Toda bracket
in the framework of stable topological model categories.
Section 6 features a comparison of different definitions of Toda brackets.  In 
Section 7, we show how universal Toda brackets are related to $k$-invariants of classifying spaces. 
The Appendix consists of a brief discussion of topological model categories and provides a technical result needed 
in Section \ref{constr_section}.

\subsection*{Acknowledgments} 
This paper is the revised version of my Ph.D. thesis at the University of Bonn \cite{Sagave_S_2006}.
I would like to thank my adviser Stefan Schwede, who suggested this project, for 
his continuous encouragement and his help in various questions arising along the way.
I also thank Christian Ausoni, Kristian Br{\"u}ning, and G\'erald Gaudens for a lot of discussions on this project. 
Moreover, I learned a lot from conversations with Teimuraz Pirashvili, and 
I benefited from discussions with Mamuka Jibladze, Henning Krause, 
Birgit Richter, John Rognes, and Brooke Shipley.

During my time as a Ph.D. student, I was supported by the Mathematical Institute of the University of Bonn,
the SFB 478 at the University of M{\"u}nster, the 
Institut Mittag-Leffler in Djursholm, a ``Kurz\-sti\-pen\-di\-um'' of the German Academic Exchange Service, 
and the GRK 1150 ``Homotopy and Cohomology'' in Bonn.

\section{Mac Lane cohomology}
We review the definition of Mac Lane cohomology via cohomology of categories and introduce a cup pairing
between the Mac Lane cohomology of a graded ring $\Lambda$ and a group of $\Lambda$-module homomorphisms. 
As a reference for Mac Lane cohomology, we recommend the last chapter of Loday's book \cite[Chapter 13]{Loday_cyclic_98}. 

\subsection{Cohomology of categories and $\HML$}
Let $\cC$ be a small category. A {\em $\cC$-bimodule} is a functor $D \colon \cC^{\op} \times \cC \to \Ab$. For
a map $f \colon X \to Y$ in $\cC$, we denote the abelian group $D(X,Y)$ by $D_f$. For maps $g \colon  X' \to X$, $h \colon  Y \to Y'$,  
and $f \colon  X \to Y$, the 
$\cC$-bimodule structure induces actions $g^*  \colon  D_f \to D_{fg}$
and $h_*  \colon  D_f \to D_{hf}$. If $A$ is a ring and $\cC$ is the category of $A$-modules, the bifunctor
$\Hom_A(-,-)$ provides an example for a $\cC$-bimodule. 

To define the cohomology a category $\cC$ with coefficients in a $\cC$-bimodule $D$ we consider the 
cochain complex $C^*(\cC, D)$ with
\[  
C^n(\cC,D) = \begin{cases} \{ c \colon N_n(\cC) \to \coprod_{g \in {\rm Mor}(\cC)} \! D_g
\; | \; c(g_1, \dots,
g_n) \in D_{g_1  \cdots g_n} \} & \text{ for $n \geq 1$}  \\
\{ c \colon \Objec (\cC) \to \coprod_{X \in {\Objec}(\cC)} \!
D(X,X) \;| \; c(X) \in D(X,X) \} & \text{ for $n=0$.} \end{cases} \]
Here $N(\cC)$ is the nerve of $\cC$, so an element of $N_n(\cC)$ is a sequence $(g_1, \dots, g_n)$ 
of $n$ composable maps in $\cC$. The abelian group structure on $C^n(\cC, D)$ is given by the pointwise addition in 
$D_g$. For $n>1$, the differential $\delta \colon  C^{n-1}(\cC, D) \to C^n(\cC, D)$ is
\[ \begin{split} (\delta c)(g_1, \dots , g_n) = & (g_1)_* c(g_2, \dots ,g_n) 
 + \sum_{i=1}^{n} (-1)^i c(g_1, \dots , g_i g_{i+1}, \dots ,g_n) \\  
 & + (-1)^{n+1} (g_n)^* c(g_1, \dots ,g_{n-1}). \end{split} \]
For $n=1$, it is $(\delta c)(g_1\colon X_1 \to X_0) = (g_1)_* c(X_1) - (g_1)^*c(X_0)$. It is easy to verify $\delta^2 = 0$.
 
\begin{definition} \cite[Definition 1.4]{Baues_W_85} 
The {\em cohomology} $H^{*}(\cC, D)$ {\em of the category} $\cC$ with coefficients in the $\cC$-bimodule $D$ 
is the cohomology of $(C^{*}(\cC, D), \delta)$. 
\end{definition}

There is a normalized version of this. A category is {\em pointed} if it has a 
preferred zero object $*$, i.e., $*$ is both initial and terminal. A {\em zero morphism} in a pointed 
category is a map which factors through the zero object. 
If $\cC$ is a pointed category, a $\cC$-bimodule $D$ is {\em normalized} if $D(*,X)=0=D(X,*)$ holds 
for all objects $X$. 

For a pointed category $\cC$ and a normalized $\cC$-bimodule $D$, we consider the subgroup 
$ \ovl{C}^n(\cC,D) = \{ c \in C^n(\cC, D) | c(g_1, \dots, g_n) = 0 \text{ if $g_i$ is zero for some $i$} \} $
of normalized cochains in $C^n(\cC, D)$. As $D$ is normalized, $\ovl{C}^*(\cC, D)$ is a subcomplex of $C^*(\cC, D)$. 
By \cite[Theorem 1.1]{Baues_D_89}, 
the inclusion $\ovl{C}^*(\cC, D) \to C^*(\cC, D)$ induces an isomorphism in cohomology. 
Therefore, we can assume representing cochains to be normalized as soon as we consider the cohomology of a {\em pointed} 
category with coefficients in a {\em normalized} bimodule.

Cohomology of categories has good naturality properties. If $F \colon \cC \to \cD$ is a functor and $D$ a
$\cD$-bimodule, there is an induced $\cC$-bimodule $F^*D$, and $F$ induces maps 
$F^* \colon C^*(\cD, D) \to C^*(\cC, F^* D)$ and $F^* \colon H^*(\cD, D) \to H^*(\cC, F^*D)$. 
The latter map is an isomorphism if $F$ is an equivalence of categories \cite[Theorem 1.11]{Baues_W_85}.

For a ring $A$, we denote the category of finitely generated free right $A$-modules by $F(A)$. 
To avoid set theoretic problems, we assume $F(A)$ to be small,
i.e., we require it to contain only one element from each isomorphism class of objects. The category 
$F(A)$ is pointed by the trivial module, 
and for an $A$-bimodule $M$, the functor $\Hom_A(-, - \tensor_A M)$ is a normalized
$F(A)$-module. 

\begin{definition}
Let $A$ be a ring and let $M$ be an $A$-bimodule. The {\em Mac Lane cohomology} of $A$ with coefficients
in $M$ is defined by 
\[ \HML^s(A, M) = H^s(F(A), \Hom_A(-, - \tensor_A M)). \]
If $M$ equals $A$, we adopt the convention $\HML^s(A) = \HML^s(A,A)$. 
\end{definition}

Mac Lane cohomology was originally defined by Mac Lane in 1956 \cite{MacLane_1956}.  Jibladze and Pirashvili \cite{Jibladze_P_1991} proved the equivalence of Mac Lane's definition to the one we use. Mac Lane cohomology is also
isomorphic to $\Ext$-groups in the abelian category  $\cF(A)$ of functors from $F(A)$ to $\modu{A}$ (see \cite{Jibladze_P_1991}) 
and to {\em topological Hochschild cohomology} (see \cite{Boekstedt_1985} for the definition and  \cite{Pirashvili_W_1992} or \cite[Theorem 6.7]{Schwede_2001} for the equivalence).
The computation of Mac Lane cohomology is known for many examples, including the cases $\HML^*(\mF_p)$ 
(see \cite{Franjou_L_S_1994}) and $\HML^*(\mZ)$ (see \cite{Franjou_P_1998}) we encounter in Section \ref{ring_sp_section}. 

For later use we prove
\begin{lemma} \label{discard_constant}
Let $A$ be a ring, let $M$ and $P$ be $A$-modules with $P$ projective, let $I \colon F(A) \to \modu{A}$ be the inclusion functor and let $T \colon F(A) \to \modu{A}$ be any functor. For $i \geq 1$, there is an isomorphism 
\[ H^i(F(A), \Hom_A(I(-) \oplus P, T(-) \oplus M)) \iso H^i(F(A), \Hom_A(-,T(-))). \]
\end{lemma}
\begin{proof}
By \cite[Corollary 3.11]{Jibladze_P_1991}, this translates to a statement about the $\Ext$-group
$\Ext_{\cF(A)}^i(I (-) \oplus P, T(-) \oplus M)$. Since the constant functor represented by $P$ is projective
in $\cF(A)$, it cancels out in the first variable. As $I$ is reduced, i.e., $I(0)=0$, it has a projective resolution by
reduced functors. Since there is only the trivial map from a reduced functor to the constant functor $M$, the 
$M$ cancels out as well. 
\end{proof}

\subsection{Mac Lane cohomology of graded rings}
If $\Lambda$ is a graded ring, the morphisms between graded $\Lambda$-modules $M$ and $N$ form a 
graded abelian group by setting $ \Hom^{i}_{\Lambda}(M,N) = \Hom_{\Lambda}(M,N[i]) = \Hom_{\Lambda}(M,N)_{-i}.$

\begin{definition} \label{split_and_sparse_def}
A graded ring, a graded abelian group,
or a graded module is {\em $n$-sparse} if it is concentrated in degrees divisible by $n$. 
A full subcategory $\cC$ of $\modu{\Lambda}$ is {\em $n$-split} if for each pair of objects $M$ and $N$ in 
$\cC$, the graded abelian group $\Hom_{\Lambda}(M,N)_*$ is $n$-sparse. 
\end{definition}

For a graded ring $\Lambda$, let 
$F(\Lambda)$ be the category of finitely generated free graded right $\Lambda$-modules. The objects of $F(\Lambda)$ 
are finite sums of shifted copies of the free module of rank $1$. 
If $\Lambda$ is $n$-sparse for $n \geq 1$,  the full subcategory of $F(\Lambda)$ given by the $n$-sparse $\Lambda$-modules
is denoted by $F(\Lambda, n)$.
For $n=1$, we have $F(\Lambda) = F(\Lambda,n)$. The category $F(\Lambda,n)$
is an example of an $n$-split subcategory of $\modu{\Lambda}$. 

\begin{definition} \label{HML_gr_def}
Let $\Lambda$ be an $n$-sparse graded ring, and let
$M$ be a graded right $\Lambda$-module. The {\em graded $n$-split Mac Lane cohomology} of $\Lambda$ 
with coefficients in $M$ is defined by
\[ \HML^{s}_{n\spl}(\Lambda,M) = H^s(F(\Lambda,n), \Hom_{\Lambda}(-, - \tensor_{\Lambda}M)). \]
If $M=\Lambda[t]$, a $t$-fold shift of $\Lambda$ for some $t \in \mZ$, we adopt the convention
\[ \HML^{s,t}_{n\spl}(\Lambda) = \HML^{s}_{n\spl}(\Lambda, \Lambda[t]). \]
If $n=1$, we drop `$1\spl$' from the notation and write $\HML^s(\Lambda, M)$ or 
$\HML^{s,t}(\Lambda)$. 
\end{definition}
The graded Mac Lane cohomology is related to the ungraded theory. If $\Lambda$ is $n$-sparse, the functor 
$- \tensor_{\Lambda_0}\Lambda \colon F(\Lambda_0) \to F(\Lambda,n)$ satisfies 
\[ (- \tensor_{\Lambda_0}\Lambda)^* \Hom_{\Lambda}(-,- \tensor_{\Lambda}\Lambda[-n]) 
\iso \Hom_{\Lambda_0}(-,-\tensor_{\Lambda_0}\Lambda_n)\]
and therefore induces a restriction map $\HML^{*,-n}_{n\spl}(\Lambda) \to \HML^*(\Lambda_0, \Lambda_n).$

A central unit $u$ of degree $n$ in $\Lambda$ is a homogeneous element $u$ of degree $n$ which is a unit and is central
in the graded sense. If $\Lambda$ has a central unit, $- \tensor_{\Lambda_0}\Lambda$ is an equivalence of categories,
and $\Lambda_n$ is isomorphic to $\Lambda_0$ as $\Lambda_0$-bimodules. This proves

\begin{lemma} \label{HML_gr_ugr}
Let $\Lambda$ be an $n$-sparse graded ring with a central unit $u$ of degree $n$. Then the restriction  
induces  an isomorphism 
$\HML^{*,-n}_{n\spl}(\Lambda) \to \HML^*(\Lambda_0)$. 
\end{lemma}

\subsection{Relation to group cohomology} We review some well known maps from Mac Lane cohomology to group cohomology. 

For an object $X$ in a category $\cC$, we denote its group of automorphisms by $\Aut(X)$.
The category with a single object $X$ and $\Hom(X,X) = \Aut(X)$ is denoted by $\cAut(X)$. It comes
with a canonical inclusion functor $\cAut(X) \to \cC$. 
If $D$ is an $\cAut(X)$-bimodule, the automorphism group $\Aut(X)$ acts via the conjugation 
action $ g x = (g^{-1})^* (g_*(x))$
from the left on the abelian group $D(X,X)$. 

\begin{proposition} \label{res_map}
Let $\cC$ be a small category, let $X$ be an object of $\cC$, and let $D$ be a $\cC$-bimodule. The 
inclusion functor $F \colon \cAut(X) \to \cC$ induces a restriction map 
\[ \Theta \colon H^*(\cC, D) \to H^*(\cAut(X),F^*D) \xrightarrow{\iso} H^*(\Aut(X), D(X,X))\]
from the cohomology of $\cC$ with coefficients in $D$ to the cohomology of the group $\Aut(X)$ with 
coefficients in the $\Aut(X)$-module $D(X,X)$. 
\end{proposition}
\begin{proof}
The first map is the restriction along the inclusion. The second map is analogous to the Mac Lane isomorphism
between the Hochschild homology of a group ring and group homology \cite[Proposition 7.4.2]{Loday_cyclic_98}.
On a cochain $c$, the isomorphism is given by 
$(\varphi(c))(g_1, \dots, g_n) = (g_n^{-1} \dotsm g_1^{-1})^* c(g_1, \dots, g_n).$
\end{proof}
When $A$ is a ring and $M$ is an $A$-bimodule, we write as usual $\GL_q A$ for the group of invertible $(q \times q)$-matrices,
which acts on the abelian group $\Mat_q M$ of all $(q \times q)$-matrices with entries in $M$ by conjugation.  
The last proposition specializes to Mac Lane cohomology for graded and ungraded rings:

\begin{corollary} \label{res_map_HML}
Let $\Lambda$ be an $n$-sparse graded ring, let $A$ be a ring, and let $M$ be an $A$-bimodule. For $q \geq 1$, 
there are restriction maps 
\[ \HML^{*,-n}_{n\spl}(\Lambda) \to H^*(\GL_q \Lambda_0, \Mat_q \Lambda_n)  \textrm{ and } 
\HML^{*}(A,M) \to H^*(\GL_q A, \Mat_q M).\]
If $A = \Lambda_0$ and $M= \Lambda_n$, the first map factors through the second map and the restriction 
$\HML^{*,-n}_{n\spl}(\Lambda) \to \HML^*(\Lambda_0, \Lambda_n).$ 
\end{corollary}

\subsection{The Cup-product}
In the following, $\Ext$-groups are understood in the sense of 
Yoneda. For a graded ring $\Lambda$, shifting of modules gives 
rise to a bigrading on $\Ext$, that is, $\Ext^{s,t}(M,N) = \Ext^{s}(M,N[t])$. 

\begin{construction} \label{cup_map}
Let $\Lambda$ be an $n$-sparse graded ring. Let
$M$ and $N$ be $\Lambda$-modules such that $M$ admits a resolution by objects in $F(\Lambda, n)$. Then there
is a well defined map 
\[\Hom_{\Lambda} (M,N) \times \HML^{s,t}_{n\spl}(\Lambda) \to 
\Ext^{s,t}_{\Lambda}(M,N), \quad (f, \gamma) \mapsto f \cup \gamma  \]
which we refer to as the {\em cup product}. It is bilinear and natural in the sense that 
$ (gf) \cup\gamma  = g_* (f \cup \gamma)$ holds for composable maps of $\Lambda$-modules $f$ and $g$. 

To define the cup product, we choose a resolution $ \dots \to M_1 \xrightarrow{\lambda_1} M_0 \xrightarrow{\lambda_0} M$
of $M$ by objects $M_i$ of $F(\Lambda,n)$ and a  normalized cocycle
\[ c \in  \ovl{C}^s(F(\Lambda,n), \Hom_{\Lambda} (-, - \otimes_{\Lambda} \Lambda[t])) \] 
representing the cohomology class $ \gamma \in \HML^{s,t}_{n\spl}(\Lambda)$.

Since $ \delta (c) = 0$ and the $\lambda_i$ form a resolution, evaluating $\delta(c)$ on $(\lambda_1, \dots ,\lambda_{s+1})$
yields 
$ \lambda_0[t]  c(\lambda_1, \dots, \lambda_s)  \lambda_{s+1}
= (-1)^s \lambda_0[t] \lambda_1[t] c(\lambda_2, \dots,
\lambda_{s+1})
= 0.$
This implies that there is a dotted arrow $\tau$ such that diagram of Figure \ref{cup_diagram_fig} commutes. 
\begin{figure}
\[ \xymatrix@-1pc{
M_{s+1} \ar[rr]^{\lambda_{s+1}} & & 
M_s \ar[dd]_{c(\lambda_1, \dots, \lambda_s)} \ar[drr] \ar[drrr]^{\lambda_s} & & & & & &\\
& & & 0 \ar[r] & \ker \lambda_{s-1} \ar@{-->}[ddll]^{\tau}  \ar[drr] \ar[r] & M_{s-1} \ar[r]^{\lambda_{s-1}} & \dots \ar[r]^{\lambda_1} & M_0 \ar[r]^{\lambda_0} & M \ar[r] & 0. \\
& & M_0[t] \ar[d]_{\lambda_{0}[t]} & & & & 0 & & \\
& & M[t]}
\] 
\caption{Defining the cup product} \label{cup_diagram_fig}
\end{figure}

If $\Psi \in \Ext_{\Lambda}^{s,0}(M, \ker \lambda_{s-1})$ denotes the Yoneda class of the
extension \[ 0 \to \ker \lambda_{s-1} \to M_{s-1} \to \dots \to M_0 \to M \to 0,\]
we define $f \cup \gamma $ to be $(-1)^{\frac{n(n+3)}{2}}((f[t])  \tau)_* (\Psi) \in \Ext_{\Lambda}^{s,t}(M,N)$.
The mysterious sign is built in to cancel out with another sign arising in Lemma \ref{filtered_Postnikov}. 
(This will keep signs out of the statements of the main results.)

The bilinearity and the naturality with respect to composition of maps are
obvious.
In Lemma \ref{independent_of_cycle} and Lemma \ref{independent_of_resolution} we show that the 
Ext-class of $f \cup \gamma $ doesn't depend on the choice of the
cocycle representing $\gamma$ and the resolution of $M$. 
\end{construction}

\begin{remark}
If $E$ is a graded $k$-algebra over a field $k$, the tensor product of a right module with a bimodule has
a left derived functor  
\[ - \tensor^L - \colon \Hom_E(P,Q) \times \HH^{s,t}_k(E) \to \Ext_E^{s,t}(P,Q). \]
Our cup-product should be thought of as similar to this. The relation becomes clearer when $\HML^*$
is defined via $\Ext$-groups in the category $\cF(A)$ of functors $F(A) \to \modu{A}$. We sketch the
ungraded case. 

The self-extensions $\Ext^*_{\cF(A)}(I,I)$ in $\cF(A)$ of the inclusion functor 
$I$ are isomorphic to $\HML^*(A)$ \cite{Jibladze_P_1991}. We can enlarge
$F(A)$ by a bigger small additive subcategory $\cC$ of $\modu{A}$ that contains $M$ without changing
the $\Ext$-group \cite[\S 2 and Corollary 3.11]{Jibladze_P_1991}. 
Evaluating an element of $\Ext_{\cC}(I,I)$ on $M$ gives an element of $\Ext^*_A(M,M)$, and 
inspecting the proof of the isomorphism $\Ext_{\cF(A)}^*(I,I) \iso \HML^*(A)$ 
\cite[Theorem B]{Jibladze_P_1991}, we see that this recovers the cup-product. We do not go into the
details as we only use the description of the product given above. 
\end{remark}

\begin{lemma} \label{trivial_maps_on_ext}
Let $ 0 \to  M' \xrightarrow{g} M_{n-1}  \to \dots \to M_0 \to M $ be an exact sequence in $\modu{\Lambda}$ 
representing a class $\Psi \in \Ext^n_{\Lambda}(M,M')$. Assume that 
$M_0, \dots ,M_{n-1}$ are free.
For maps $f \colon M' \to N$ and $h \colon M_{n-1} \to N$, we have 
$(f + h g)_* (\Psi) = f_* (\Psi)$.
\end{lemma}
\begin{proof}
This statement becomes trivial with $\Ext$ defined via projective resolutions.
\end{proof}

\begin{lemma} \label{independent_of_cycle}
The cup product of Construction \ref{cup_map} does not depend on the choice of the cocycle
representing $\gamma$. 
\end{lemma}
\begin{proof}
It is enough to show that the extension associated to a coboundary 
represents the trivial element in $\Ext_{\Lambda}^{s,t}(M,N)$.
Let $b \in  \ovl{C}^{s-1}(F(\Lambda,n), \Hom_{\Lambda} (-, - \otimes_{\Lambda} \Lambda[t]))$ be a
normalized cochain. Evaluating $\delta(b)$ on $(\lambda_1, \dots, \lambda_s)$ yields 
$\lambda_0 [t]  \delta(b)( \lambda_1, \dots, \lambda_s) = 
(-1)^{s} \lambda_0 [t] b( \lambda_1, \dots, \lambda_{s-1}) \lambda_s.$
Hence the $\tau$ associated to $c=\delta(b)$ extends to $M_{s-1}$, so $((f[t]) \tau)_* (\Psi) = 0$ 
by the last lemma.  
\end{proof}

\begin{lemma} \label{independent_of_resolution}
The cup product of Construction \ref{cup_map} does not depend on the choice of the
resolution of $M$. 
\end{lemma}
\begin{proof}
Suppose we are given another resolution 
$ \dots \xrightarrow{\lambda'_1}  M'_0 \xrightarrow{\lambda'_0} M$
of $M$ by objects of $F(\Lambda, n)$. 
Then there exist maps $\alpha_i \colon M'_i \to M_i$ with $\lambda_i \alpha_i = \alpha_{i-1} \lambda'_i$ and
$\lambda_0 \alpha_0 = \lambda'_0$. The problem is that in general 
$ (\lambda'_0[t]) (c(\lambda'_1, \dots, \lambda'_s))=  
(\lambda_0[t])(c(\lambda_1, \dots, \lambda_s)) \alpha_s$ does not hold. 
As we are only interested in the induced
maps on Ext-groups, it suffices to show that the two maps  
give rise to maps $\tau, \tau' \colon \ker \lambda'_{s-1} \to M[t]$ which induce the same map 
$\Ext_{\Lambda}^{s,0}(M, \ker \lambda'_{s-1}) \to \Ext_{\Lambda}^{s,t}(M,N).$ 

Using $\lambda_i \alpha_i = \alpha_{i-1} \lambda'_i$ and the definition of $\delta$, we
obtain the equation 
\begin{align*}
 0 &=  (\delta c)(\alpha_0 , \lambda'_1, \dots, \lambda'_s)
      + (-1)^{s} (\delta c)(\lambda_1, \dots, \lambda_s, \alpha_s) \\
   & \quad + \sum_{i=1}^{s-1} (-1)^i (\delta c)(\lambda_1, \dots, 
     \lambda_i, \alpha_i, \lambda'_{i+1}, \dots, \lambda'_{s}) \\
   &= \alpha_0[t]  c(\lambda'_1, \dots, \lambda'_s) -
     c(\lambda_1, \dots, \lambda_s)  \alpha_s + \lambda_1[t]  g + h  \lambda'_s
\end{align*}
in which  $h\colon M'_{s-1} \to M_0[t]$ and $g\colon M'_s \to M_1[t]$ are maps 
we don't need to know explicitly. Composing with $\lambda_0[t]$  and applying Lemma \ref{trivial_maps_on_ext}
completes the proof. 
\end{proof}

\section{Realizability in triangulated categories}
In this section we give a quick review of the obstruction theory for realizability in triangulated categories described in 
\cite[Appendix A]{Benson_K_S_04}. The necessary background on triangulated categories can be 
found in  Weibel's book \cite{Weibel-intorduction_1994}.

Let $\cT$ be a triangulated category, which we always assume to have infinite coproducts. An object 
$N$ of $\cT$ is {\em compact} if the functor $\cT(N,-)$ preserves arbitrary coproducts. For objects
$X$ and $Y$ of $\cT$, we write $\cT(X,Y)_*$ for the graded abelian group whose degree $k$ part is
$\cT(X[k],Y)$. 

We fix a compact object $N$ in $\cT$. Under composition, $\Lambda := \cT(N,N)_*$ becomes a graded ring,
and $\cT(N,X)_*$ is a right $\Lambda$-module for 
every object $X$. The resulting functor $\cT(N,-) \colon \cT \to \modu{\Lambda}$ from 
$\cT$ to graded $\Lambda$-modules maps distinguished triangles in $\cT$ to long exact sequences. 
Furthermore, it preserves arbitrary coproducts since $N$ is compact, and it commutes with
the shift of $\cT$ and $\modu{\Lambda}$. 

\begin{definition}
In the above context, a $\Lambda$-module $M$ is called {\em realizable} if there is an object
$X$ in $\cT$ such that $\cT(N,X)_* \iso M$.
\end{definition}

The following example for this situation is studied in \cite{Benson_K_S_04}. Let $A$ be a differential graded algebra 
over a field $k$, and let $\cT = \cD(A)$  be the derived category of dg $A$-modules. If $N$ is the free
module of rank $1$, we have $\cT(N,N)_{-*} = H^*(A)$. The realizability question amounts
to whether a graded module over the cohomology ring $H^*(A)$ is the cohomology of a dg
$A$-module. In Section \ref{ring_sp_section}, we will address the corresponding question for a ring spectrum $R$:
when is a module over the homotopy groups of $R$ the homotopy of an $R$-module spectrum?

\subsection{Realizability obstructions}
Let $\cT$, $N$ and $\Lambda$ be as above. An object of $\cT$ is called {\em $N$-free} 
if it is a sum of shifted copies of $N$. We note that $\cT(N,-)_*$ restricts to an equivalence between the full 
subcategory of $N$-free objects in $\cT$ and the category of free $\Lambda$-modules. 

\begin{definition} \label{postnikov_sys_def} \cite[Definition A.6]{Benson_K_S_04}
For $k \geq 1$, an {\em $N$-exact $k$-Postnikov system} for a $\Lambda$-module $M$ consists of an epimorphism 
$\cT(N,X_0)_* \to M$ and a diagram 

\[\xymatrix{ & Y_{k-1} \ar[d]_(.6){\pi_{k-1}} & Y_{k-2} \ar[d]_(.6){\pi_{k-2}} \ar@{+->}[l]_{\alpha_{k-1}} 
& Y_{k-3} \ar@{+->}[l]_{\alpha_{k-2}}   & Y_2 \ar[d]_(.6){\pi_{2}} 
& Y_1 \ar[d]_(.6){\pi_{1}} \ar@{+->}[l]_{\alpha_{2}} & Y_0 = X_0 \ar@{+->}[l]_(.6){\alpha_{1}} \\
 X_k \ar[ur]^(.4){\iota_k} & X_{k-1} \ar[ur]_(.4){\iota_{k-1}} 
& X_{k-2} \ar[ur]_(.4){\iota_{k-2}} \ar@{}[urrr]|{\dots} &  & X_2  \ar[ur]_(.4){\iota_{2}} & X_1  
\ar[ur]_(.4){\iota_{1}}}\]
such that all arrows of the form $\xymatrix@1{& \ar@{+->}[l]}$ denote 
morphisms of degree $1$, all triangles are distinguished triangles in $\cT$, and
each object $X_i$ is $N$-free. Moreover, the maps $d_j = \pi_{j-1} \iota_{j}$ with $j \geq 2$ and $d_1 = \iota_1$ are required to 
induce an exact sequence
\[ \cT(N,X_k)_* \xrightarrow{(d_k)_*} \cT(N,X_{k-1}) \xrightarrow{(d_{k-1})_*} 
\dots \xrightarrow{(d_1)_*} \cT(N, X_0)_{*} \to M \to 0. \]
An {\em $N$-exact Postnikov system} is a collection of distinguished triangles as above which extends 
infinitely to the left.
\end{definition}

Proposition A.19 of \cite{Benson_K_S_04} shows that a $\Lambda$-module $M$ is realizable if there exists an 
$N$-exact Postnikov system of $M$. By realizing the first two steps of a free resolution of $M$, one 
can easily see that $N$-exact $2$-Postnikov systems exist for every $M$. Therefore, the realizability
problem can be attacked by extending Postnikov systems stepwise to the left. 

By \cite[Lemma A.12(iii)]{Benson_K_S_04}, every $N$-exact $k$-Postnikov system of $M$ induces an exact sequence
\begin{multline*}
\cT(N,X_1)_*[1-k] \xrightarrow{(d_1)_*} \cT(N,X_0)_*[1-k] \xrightarrow{\alpha_*} \cT(N, Y_{k-1})_* \\
\xrightarrow{(\pi_{k-1})_*} \cT(N, X_{k-1})_* \xrightarrow{(d_{k-1})_*} \cT(N,X_{k-2})_* 
\end{multline*}
of $\Lambda$-modules, where the map $\alpha\colon X_0[1-k] = Y_0[1-k] \to Y_{k-1}$ is the 
composition $\alpha_{k-1} \dotsm \alpha_1$. Hence there is an exact sequence 
\setcounter{equation}{\value{subsection}}
\begin{equation} \label{representing_obstr_cls} \stepcounter{subsection}
\begin{split}
 0 \to M[1-k] \xrightarrow{\eta_{k-1}} \cT(N,Y_{k-1})_* \xrightarrow{(\pi_{k-1})_*} \cT(N, X_{k-1})_* 
\xrightarrow{(d_{k-1})_*} \dots \\ \dots \xrightarrow{(d_2)_*} \cT(N,X_1)_* 
\xrightarrow{(d_1)_*} \cT(N,X_0)_* \to M \to 0. \end{split} \end{equation}
of $\Lambda$-modules. Its Yoneda class is denoted by $\kappa_{k+1}(M) \in \Ext_{\Lambda}^{k+1,k-1}(M,M)$
and is called the  obstruction class associated to the Postnikov system because of 

\begin{lemma} \cite[Lemma A.18]{Benson_K_S_04} \label{existence_postn}
If the class $\kappa_{k+1}(M)$ of an $N$-exact $k$-Postnikov system of $M$ is trivial, then there
exists an $N$-exact $(k+1)$-Postnikov system for $M$ whose underlying $(k-1)$-Postnikov system agrees with
that of the given one. 
\end{lemma}

The class $\kappa_3(M)$ is always defined and unique \cite[Proposition 3.4(ii)]{Benson_K_S_04}. 
If the higher obstructions $\kappa_i(M)$ for $i \geq 4$ are defined, they may depend on the choice of the Postnikov
system. 

\subsection{A criterion for uniqueness of obstruction classes}
To compare the obstruction classes of different Postnikov systems, we need
\begin{definition} \label{map_of_Postnikov}
Let $(X_j, Y_j, \alpha_j, \iota_j, \pi_j, M)$ and $(X'_j, Y'_j, \alpha'_j, \iota'_j, \pi'_j, M)$ be two $N$-exact
$k$-Postnikov systems for $M$. A morphism between them consists of maps $f_j\colon X_j \to X'_j$ and $g_j \colon Y_j \to Y'_j$
such that $f_{k-1} d_k = d'_k f_k$ and the following commutativity relations hold for $1 \leq j \leq k-1$:
\[
g_{j-1} \iota_j = \iota'_j f_j \qquad (g_j[1]) \alpha_j = \alpha'_j g_{j-1} \qquad 
f_j \pi_j = \pi'_j g_j.  
\]
More generally, for $1 \leq l \leq k$, an $l$-map of $N$-exact $k$-Postnikov systems for $M$ is a 
map of the underlying $N$-exact $l$-Postnikov systems. 
\end{definition}
A map between two
$N$-exact $k$-Postnikov systems induces a map of the long exact sequences representing the obstruction 
classes, and this map is $\id_M$ on the outer terms. So the obstruction classes of two Postnikov systems 
coincide if there is a map between them. 
Note that this does {\em not} need the relation 
$g_{k-1} \iota_k = \iota'_k f_k$, which therefore wasn't required in Definition \ref{map_of_Postnikov}.
To produce such maps, we use

\begin{lemma} \label{uniqueness_postn}
Suppose we are given an $l$-map between two $N$-exact $k$-Postnikov systems with $1 \leq l < k$. There
is an element in $\Ext^{l,1-l}_{\Lambda}(M,M)$ whose vanishing implies the existence of an $(l+1)$-map between
the Postnikov systems.
\end{lemma}
\begin{proof}
Since both $(\cT(N,X_i)_*, (d_i)_*)$ and $(\cT(N,X'_i)_*, (d'_i)_*)$ are exact complexes of free $\Lambda$-modules,
we can find a $f_{l+1} \colon X_{l+1} \to X'_{l+1}$ with $f_l d_{l+1} = d'_{l+1} f_{l+1}$.

Let us assume for a moment  our map of Postnikov systems satisfies 
$g_{l-1} \iota_l = \iota'_l f_l$. Then we could find a $g_l \colon Y_l \to Y'_l$ such that 
$(f_l, g_{l-1}, g_l[1])$ is a map between the distinguished triangles $( - \pi_l[1], \alpha_l, \iota_l)$ and
$( - \pi'_l[1], \alpha'_l, \iota'_l)$. The maps $f_{l+1}$ and $g_l$ would complete the required data of 
an $(l+1)$-map.

In general, $\varphi = \iota'_l f_l - g_{l-1} \iota_l \in \cT(X_l,Y'_{l-1})$ is non-zero. By applying
$\cT(X_l,-)$ to the triangle $(\pi'_l, -\alpha'_l[-1], -\iota'_l[-1])$ we see that
there is a $\psi \in \cT(X_l, Y'_{l-2}[-1])$ with $(\alpha'_{l-1})_* (\psi) = \varphi$. 

If we apply $\cT(N,-)_*$ to $\cT(X_l, Y'_{l-2}[-1])$ and use \cite[Lemma A.12(i)-(ii)]{Benson_K_S_04},  
$\psi$ defines a class in $\Ext^{l,1-l}_{\Lambda}(M,M)$. The vanishing of this $\Ext$-group
implies the existence of a $\rho \in \cT(X_{l-1}, Y'_{l-2}[-1])$ with
$\rho d_l =  \psi$.

Now we can change our map of Postnikov systems by 
replacing $g_{l-1}$ by $\ovl{g_{l-1}} = g_{l-1} + (\alpha'_{l-1}[-1]) \rho \pi_{l-1}$. The  $\ovl{g_{l-1}}$
satisfies the required relations. In addition, 
\[ 
\ovl{g_{l-1}} \iota_l
= g_{l-1} \iota_l+ (\alpha'_{l-1}[-1]) \rho d_l 
= g_{l-1} \iota_l + \varphi = \iota'_l f_l
\]
holds, and the modified $l$-map extends to an $(l+1)$-map by the argument above. \end{proof}

Recall that a graded abelian group or a graded ring is {\em $n$-sparse} if it is concentrated in degrees divisible by $n$. 

\begin{corollary} \label{uniqueness_obstr}
Suppose that $\Lambda=\cT(N,N)_*$ is $n$-sparse and  $M$ is an $n$-sparse $\Lambda$-module. Then there exists
an $N$-exact $(n+1)$-Postnikov system of $M$, and all $N$-exact $(n+1)$-Postnikov systems of $M$ give rise to
the same obstruction class $\kappa_{n+2}(M) \in \Ext_{\Lambda}^{n+2,-n}(M,M)$.   
\end{corollary}
\begin{proof}
The groups $\Ext^{l+1,1-l}_{\Lambda}(M,M)$ vanish for $2 \leq l \leq n$ because of the sparseness of $\Lambda$ and $M$. 
Hence there is an $N$-exact $(n+1)$-Postnikov system for $M$ by Lemma \ref{existence_postn}. Similarly, Lemma 
\ref{uniqueness_postn} and the
vanishing of $\Ext^{l,1-l}_{\Lambda}(M,M)$ for $2 \leq l \leq n$ provide the existence of a map between 
two $N$-exact $(n+1)$-Postnikov systems for $M$. This implies the uniqueness of the obstruction class $\kappa_{n+2}(M)$.  
\end{proof}

\section{Toda brackets and realizability}
We recall the definition of Toda brackets in triangulated categories and show how they are related 
to the realizability obstructions of the last section. 

\subsection{Definition of higher Toda brackets}
Cohen's definition \cite[\S 2]{Cohen_1968} of higher Toda brackets can be interpreted in the context
of triangulated categories. We follow Shipley \cite[Appendix A]{Shipley_2002} in doing so. 

\begin{definition} \label{filtered_obj} \cite[Definition A.1]{Shipley_2002}
Let $\cT$ be a triangulated category and let
\[X_{n-1} \xrightarrow{\lambda_{n-1}} X_{n-2} \xrightarrow{\lambda_{n-2}} \dots \xrightarrow{\lambda_1} X_0 \]
be $(n-1)$ composable maps in $\cT$.
An {\em $n$-filtered object} $X \in \{\lambda_1, \dots, \lambda_{n-1}\}$ consists of a
sequence of maps $\ast = F_0 X \xrightarrow{i_0} F_1 X \xrightarrow{i_1} \dots \xrightarrow{i_{n}} F_{n}X = X$
and
choices of distinguished triangles 
$F_j X \xrightarrow{i_j} F_{j+1} X \xrightarrow{p_{j+1}} X_j[j] \xrightarrow{d_j} (F_j X)[1] $
such that $(p_j[1])(d_j) = \lambda_j[j]$.

The maps $X_0 \iso F_1 X \to X$ and 
$X = F_{n} X \xrightarrow{p_{n}} X_{n-1}[n-1]$ are denoted by $\sigma'_X$ and 
$\sigma_X$. 
\end{definition}
Our definition differs from \cite[Definition A.1]{Shipley_2002} in that 
we require the objects $X_j[j]$ to {\em be} the cones of the maps $i_j$, rather than to be isomorphic to the cones.
This does not make a difference since triangles isomorphic to distinguished triangles are distinguished again.

For a map $\lambda_1 \colon X_1 \to X_0$ in $\cT$, the cone $C$ of $\lambda_1$ is part of a distinguished triangle 
$X_1 \to X_0 \to C \to X_1[1]$. With the filtration $* \to X_0 \to C$, it is a $2$-filtered object 
in $\{\lambda_1\}$. 

If there exists an $n$-filtered object $X \in \{\lambda_1, \dots, \lambda_{n-1}\}$, each twofold composition
$\lambda_i \lambda_{i+1}$ has to be zero since it can be written as a composition of maps which 
contains two consecutive maps in a distinguished triangle. 

Though a filtered object consists of similar data as a Postnikov system, we emphasize the difference:
a filtered object starts from a {\em fixed complex} of maps, while a Postnikov system starts from a {\em module}
and is assumed to have {\em some} underlying {\em resolution}. Lemma \ref{filtered_Postnikov} shows how in special 
cases a filtered object gives rise to a Postnikov system.

We will construct filtered objects using
\begin{lemma} \cite[Lemma A.4]{Shipley_2002} \label{cones_filtered}
Let $\lambda_i\colon X_i \to X_{i-1}$ be a  sequence of composable maps in a triangulated category $\cT$. 
An $n$-filtered object $X \in \{ \lambda_2, \dots, \lambda_n\}$ with a map $\alpha \colon X \to X_0$ 
gives rise to an $(n+1)$-filtered object $C_{\alpha} \in \{ \alpha \sigma'_X, \lambda_2, \dots, \lambda_n \}$,
and an $n$-filtered object $X \in \{ \lambda_1, \dots, \lambda_{n-1}\}$ with a map $\alpha \colon X_n[n-1] \to X$
gives rise to an $(n+1)$-filtered object $C_{\alpha} \in \{ \lambda_1, \dots, \lambda_{n-1}, (\sigma_X \alpha)[-n+1]\}$.
\end{lemma}
\begin{proof}
The first part uses the octahedral axiom. The second part is immediate. 
\end{proof}

\begin{definition} \label{higher_tbracket} \cite[Definition A.2]{Shipley_2002}
Let $\cT$ be a triangulated category. A map $\gamma \in \cT(X_n[n-2], X_0)$ lies 
in the {\em $n$-fold Toda bracket} of  
\[X_{n} \xrightarrow{\lambda_{n}} X_{n-1} \xrightarrow{\lambda_{n-1}} \dots \xrightarrow{\lambda_1} X_0 \]
if there exist an $(n-1)$-filtered object $X \in \{ \lambda_2, \dots, \lambda_{n-1} \}$ and maps $\gamma_0\colon X \to X_0$ and 
$\gamma_n\colon X_{n}[n-2] \to X$ such that $\gamma = \gamma_0 \gamma_n$ holds 
and the two triangles in the following diagram commute:
\[ \xymatrix{ & X_1 \ar[d]_{\sigma'_X} \ar[rd]^{\lambda_1} & \\
X_n[n-2] \ar[r]^{\gamma_n} \ar[dr]_{\lambda_n[n-2]} & X \ar[d]^{\sigma_X} \ar[r]_{\gamma_0} & X_0 \\
& X_{n-1}[n-2] }\]
We write $\tbracket{\lambda_1, \dots,\lambda_n} \subseteq \cT(X_n[n-2],X_0)$ for the possibly empty
set of all those $\gamma$.\end{definition}
We refer to Remark \ref{discuss_different_tbrackets} for a discussion of other definitions of Toda brackets.

For $n=3$, this defines the triple Toda bracket $\tbracket{\lambda_1, \lambda_2, \lambda_3}$. 
The cone of $\lambda_2$ serves as the $2$-filtered object. The set
$\tbracket{\lambda_1, \lambda_2, \lambda_3}$ is non-empty iff $\lambda_1 \lambda_2 = 0 = \lambda_2 \lambda_3$. 
It is easy to check that two elements of 
$\tbracket{\lambda_1, \lambda_2, \lambda_3}$ differ by an element of the set
$(\lambda_1)_* (\cT(X_3[1],X_1)) + (\lambda_3[1])^* (\cT(X_2[1],X_0))$, 
which we refer to as the {\em indeterminacy} of the Toda bracket.

\begin{remark} \label{different_triple}
In the situation of $(\lambda_1, \lambda_2, \lambda_3)$ with $\lambda_1 \lambda_2 = 0 = \lambda_2 \lambda_3$, 
there are two more equivalent definitions of 
$\tbracket{\lambda_1, \lambda_2, \lambda_3}$ 
which involve distinguished triangles containing $\lambda_1$ or $\lambda_3$ instead of $\lambda_2$. By 
choosing distinguished triangles in the horizontal lines and appropriate extensions, one builds the commutative diagram
of Figure \ref{differen_triple_fig}.
\begin{figure}
\[ \xymatrix{
X_3 \ar@{-->}[d]_{\tau_3[-1]} \ar[r]^{\lambda_3} & X_2 \ar@{=}[d] \ar[r]^{\iota_3} & C_3 \ar@{-->}[d]^{\epsilon_3}
\ar[r]^{\pi_3} & X_3[1] \ar@{-->}[d]^{\tau_3} \ar[r]^{- \lambda_3[1]} & X_2[1] \ar@{=}[d] \\
C_2[-1] \ar@{-->}[d]_{\epsilon_2[-1]} \ar[r]^{- \pi_2[-1]} & X_2 \ar@{-->}[d]^{\tau_2[-1]} \ar[r]^{\lambda_2}
 & X_1 \ar@{=}[d] \ar[r]^{\iota_2} &
C_2 \ar@{-->}[d]^{\epsilon_2} \ar[r]^{\pi_2} & X_2[1] \ar@{-->}[d]^{\tau_2} \\
X_0[-1] \ar[r]^{- \iota_1 [-1]} & C_1[-1] \ar[r]^-{- \pi_1[-1]}  & X_1 \ar[r]^{\lambda_1} & X_0 \ar[r]^{\iota_1} & C_1.}\] 
\caption{Different definitions of triple Toda brackets} \label{differen_triple_fig}
\end{figure}
Considering the middle line as a filtered object, one sees that $\epsilon_2 \tau_3 \in \tbracket{\lambda_1, \lambda_2, \lambda_3}$ in the sense of Definition \ref{higher_tbracket} above. Starting with the upper line, one can 
first choose $\epsilon_3$. Since $\epsilon_2 \tau_3$ is a choice for extending $\lambda_1 \epsilon_3$ to 
$X_3[1]$, this is an equivalent definition not involving $C_2$. A third definition uses the
distinguished triangle in the lower line. 
\end{remark}

\subsection{Existence and indeterminacy of higher Toda brackets}
A sequence $(\lambda_1, \dots, \lambda_n)$ of composable maps has to satisfy restrictive conditions for 
its Toda bracket to be non-empty. For example, $0 \in \tbracket{\lambda_2, \dots, \lambda_{n-1}}$ is a 
necessary condition for the existence of an $(n-1)$-filtered object $X \in \{\lambda_2, \dots, \lambda_{n-1}\}$
\cite[Proposition A.5]{Shipley_2002}, and the additional 
requirement $\lambda_1 \lambda_2 = 0 = \lambda_{n-1} \lambda_n$
will in general not be sufficient for $\tbracket{\lambda_1, \dots, \lambda_n}$ to be non-empty.
We introduce an additional assumption to obtain non-empty Toda brackets with controllable indeterminacy. 

\begin{definition}
A full subcategory $\cU$ of a triangulated category $\cT$ is $n$-split if $\cT(X,Y)_*$ is 
$n$-sparse for all objects  $X$ and $Y$ of $\cU$. 
\end{definition}
This is the analog to Definition \ref{split_and_sparse_def} for triangulated categories.
If $\cT$ has a compact object $N$ for which $\cT(N,N)_*$ is $n$-sparse, the subcategory of sums of copies of 
$N$ which are shifted by integral multiples of $n$ is $n$-split. 

\begin{lemma} \label{mapping_into_filtered}
Let $\cU$ be an $n$-split subcategory of a triangulated category $\cT$ with $n \geq 2$, let
$ X_{l-1} \xrightarrow{\lambda_{l-1}}  \dots \xrightarrow{\lambda_1} X_0 $
be a sequence of maps in $\cU$ with $2 \leq l \leq n-1$, and let 
$X \in \{\lambda_1, \dots, \lambda_{l-1} \}$ be an $l$-filtered object. Then for every object $Y$ in $\cU$,
we have $\cT(Y[l],X) = 0$  and $\cT(X,Y[-1]) = 0$. 
\end{lemma}
\begin{proof}
To show the first part, we choose a map $\alpha\colon Y[l] \to X$. Since the composition 
$Y[l] \to X = F_{l}X \xrightarrow{\sigma_X} X_{l-1}[l-1]$ is zero, 
$\alpha$ factors through $F_{l-1} X \to F_{l} X$. 
Using inductively that $\cT(Y[l], X_j[j])=0$ for $j = {l-2}, \dots, 0$, we 
obtain that $\alpha$ factors through $F_0X \to F_{l}X$. Hence 
$\alpha = 0$ since $F_0 X = *$.

For the second part, we first observe that  $\cT(F_1 X, Y[-1]) \iso \cT(X_0, Y[-1]) = 0$. The exact sequence 
$ \cT(X_j[j], Y[-1]) \to \cT(F_{j+1}X, Y[-1]) \to \cT(F_j X, Y[-1])$
in which the first term is trivial for $j \leq l-2$ can be used to show the assertion by induction. 
\end{proof}

\begin{lemma} \label{exis_filtered}
Let $\cU$ be an $n$-split subcategory of a triangulated category $\cT$. Then a sequence 
$X_l \xrightarrow{\lambda_l} X_{l-1}  \xrightarrow{\lambda_{l-1}} \dots \xrightarrow{\lambda_1} X_0$
in $\cU$ with $\lambda_i \lambda_{i+1} = 0$ 
admits an $(l+1)$-filtered object $X \in \{\lambda_1, \dots ,\lambda_l \}$ if $l \leq n+1$. 
If $l \leq n$, the $(l+1)$-filtered object is unique up to isomorphism.
\end{lemma}
\begin{proof}
The map from $X_0$ to the cone of $\lambda_1\colon X_1 \to X_0$ gives the data of 
a $2$-filtered object in $\{ \lambda_1 \}$. Inductively, we 
assume that $X \in \{\lambda_1, \dots ,\lambda_{j-1} \}$ is a $j$-filtered object with $j \leq n$ and consider the solid arrow diagram
\[\xymatrix@-0.5pc{ F_{j-2} X \ar[rr]^{i_{j-1}} & & F_{j-1} X \ar[dl]^{p_{j-1}} \ar[rr]^{i_{j}} & & F_j X \ar[dl]^{p_j}& \\
&  X_{j-2}[j-2]  \ar@{+->}[ul]^{d_{j-2}}  & &  X_{j-1}[j-1] \ar@{+->}[ul]^{d_{j-1}} \ar@{+->}[ll]^{\lambda_{j-1}[j-1]}
& & \ar[ll]^{\lambda_j[j-1]} \ar@{-->}[ul]_{\beta} X_j[j-1] .} \]   
The map $(p_{j-1} d_{j-1}) (\lambda_j[j-1])$ is trivial as a shift of $\lambda_{j-1} \lambda_j$. 
Hence $d_{j-1} (\lambda_{j}[j-1])$ lifts along $i_{j-1}$ and factors through $F_{j-2}X$. We 
have $\cT(X_{j-1} [j-1], F_{j-2}X) = 0$ by 
the last lemma, hence  $d_{j-1} (\lambda_{j}[j-1]) = 0$. This provides the existence of the dotted arrow $\beta$.
By Lemma \ref{cones_filtered}, the cone of $\beta$ is a  $(j+1)$-filtered object in 
$\{\lambda_1, \dots, \lambda_j\}$. 

We prove uniqueness by inductively constructing isomorphisms $f_j \colon F_jX \to F'_jX$ compatible
with all structure maps. This is trivial for the $1$-filtered objects. Assume we are given an
 isomorphism $f_{j-1} \colon F_{j-1}X \to F'_{j-1}X$. The compatibility yields 
$p'_{j-1} d'_{j-1} = \lambda_{j-1}[j-1] = p'_{j-1} f_{j-1} d_{j-1}$. Hence the exact sequence 
resulting from 
applying $\cT(X_{j-1}[j-2], -)$ to \[F'_{j-2}X \xrightarrow{\iota'_{j-1}} F'_{j-1}X \xrightarrow{p'_{j-1}} X_{j-2}[j-2]\]
shows that $d'_{j-1}-f_{j-1} d_{j-1}$ is in the image of $(\iota'_{j-1})_*$. Since 
$\cT(X_{j-1}[j-2], F_{j-2} X')$ is trivial for $j \leq n+1$ by Lemma \ref{mapping_into_filtered}, this
implies $d'_{j-1}=f_{j-1} d_{j-1}$. Completing $(\id_{X_{j-1}[j-1]}, f_{j-1})$ to a map of triangles
yields the desired $f_j$. 
\end{proof}

\begin{proposition} \label{indet+existence}
Let $\cU$ be an $n$-split subcategory of a triangulated category $\cT$ and let 
\[ X_{n+2} \xrightarrow{\lambda_{n+2}} X_{n+1} \xrightarrow{\lambda_{n+1}}  \dots \xrightarrow{\lambda_1} X_0 \]
be a sequence of maps in $\cU$ with $\lambda_{i} \lambda_{i+1} = 0$. 
Then the Toda bracket $\tbracket{\lambda_1, \dots, \lambda_{n+2}}$ is defined, is non-empty,
and has the indeterminacy \[(\lambda_1)_{*} (\cT(X_{n+2}[n], X_1)) + (\lambda_{n+2}[n])^* (\cT(X_{n+1}[n],X_0)).\] 
\end{proposition}
\begin{proof}
An  $(n+1)$-filtered object $X \in \{\lambda_2 , \dots , \lambda_{n+1}\}$
 exists and is unique by Lemma \ref{exis_filtered}. 
To construct $\gamma_{n+2}$, we consider the  exact sequence 
\[
\cT(X_{n+1}[n], X) \xrightarrow{(\sigma_X)_*} \cT(X_{n+2}[n], X_{n+1}[n]) \to \cT(X_{n+1}[n], F_nX [1]).  
\]
The last term is trivial by Lemma \ref{mapping_into_filtered}. Hence there is a $\gamma_{n+2}$ with 
$\sigma_X \gamma_{n+2} = \lambda_{n+2}[n]$. 

To obtain $\gamma_0$, we use  $F_1 X \xrightarrow{\iso} X_1$ and $\lambda_1$ to get a map $F_1 X \to X_0$. It
can be extended to $F_2X$ since $\lambda_1 \lambda_2 = 0$. Inductively, we can extend
it to a map $\gamma_0 \colon X = F_{n+1} X \to X_0$: the obstruction for extending a map $F_{j-1}X \to X_0$
to $F_{j} X$ lies in $\cT(X_{j-1}[j-2],X_0)$, which is
trivial for $3 \leq j \leq n+1$.   

Next we compute the indeterminacy. Since we have an exact sequence 
\[\cT(X_{n+2}[n],F_{n}X) \xrightarrow{(i_n)_*} \cT(X_{n+2}[n], F_{n+1}X) \xrightarrow{(\sigma_X)_*} 
\cT(X_{n+2}[n],X_{n+1}[n]),\]
we know that two different choices of $\gamma_{n+2}$ differ by an element in the image of $(i_n)_*$.
Using the same argument as in Lemma \ref{mapping_into_filtered}, we see that every map $X_{n+2}[n] \to F_{n}X$
factors through $\sigma'_X\colon X_1 \iso F_1X \to F_n X$. Therefore, the possible difference is in the image
of $(\sigma'_X)_*$, and after composing with any choice for $\gamma_0$ we obtain that this part of the 
indeterminacy is $(\lambda_1)_{*} \cT(X_{n+2}[n], X_1)$. 

To examine the other part of the indeterminacy, we first construct an auxiliary $n$-filtered object 
$F'_nX \in \{\lambda_{3}[1], \dots , \lambda_{n+1}[1]\}$. For $0 \leq j \leq n$, we define $F'_jX$ to be part of 
a distinguished triangle 
$X_1 \to F_{j+1}X \to F'_j X$.

The exact sequence
\[\cT(F'_n X , X_0) \to \cT(X, X_0) \to \cT(X_1, X_0) \]
shows that the difference $\ovl{\gamma}_0$ of two  choices for $\gamma_0$ is in the image of $\cT(F'_{n}X, X_0)$. 
Since $\cT(F'_{n-1}X, X_0)$ vanishes by Lemma \ref{mapping_into_filtered}, 
there is an $\omega \colon X_{n+1}[n] \to X_0$
with $\omega \sigma_X = \ovl{\gamma}_0$. If we apply $(\gamma_{n+2})^*$ to $\omega \sigma_X$, we see
that this part of the indeterminacy is given by $(\lambda_{n+2})^*(\cT(X_{n+1}[n], X_0))$.
\end{proof}

\subsection{Relation to realizability obstructions}
In this section we exhibit the link between Toda brackets and realizability obstructions. More precisely,
we use the cup product of Construction \ref{cup_map} to turn the slogan `the Toda brackets of the resolution are 
realizability obstructions' into a theorem. The first step is the relation between filtered objects in the
sense of Definition \ref{filtered_obj} and Postnikov systems in the sense of Definition \ref{postnikov_sys_def}. 
 
\begin{lemma} \label{filtered_Postnikov}
Let $X_n \xrightarrow{\lambda_n} X_{n-1} \xrightarrow{\lambda_{n-1}} \dots \xrightarrow{\lambda_1} X_0$ be
a sequence of maps in $\cT$ such that each $X_i$ is $N$-free and $\cT(N,-)$ maps it to an exact sequence
of $\Lambda$-modules. Let $M$ be the cokernel of the map $(\lambda_1)_* \colon \cT(N,X_1)_* \to \cT(N,X_0)$. 
An $(n+1)$-filtered object $X \in \{\lambda_1, \dots, \lambda_n\}$ determines all data of an $N$-exact
$(n+1)$-Postnikov system of $M$ except the map $X_{n+1} \to Y_n$ by 
setting $Y_l = (F_{l+1})[-l]$ for $0 \leq l \leq n$ and 
\[ \pi_l =(-1)^l p_{l+1}[-l], \qquad \iota_l = (-1)^l d_l[-l], \qquad \textrm{ and } \qquad \alpha_l = (-1)^{l+1} i_l[-l+1] \]
for $1 \leq l \leq n$. Then $\alpha = (-1)^{\frac{n(n+3)}{2}} \sigma'_{X}[-n] $, and the $(-1)(\lambda_i)_*$ form the underlying resolution
of the Postnikov system. 
\end{lemma}
\begin{proof}
The triangles $(\alpha_l, \iota_l, \pi_l)$ are distinguished since the $(d_l, p_{l+1}, i_l)$ are. The signs needed for this 
imply $\pi_{l-1}\iota_l =  - \lambda_l$ as well as the sign relating $\alpha = \alpha_n \dotsm \alpha_1$ 
to $\sigma'_X = i_{n} \dotsm i_1$. 
\end{proof}

Before stating the main theorem of this section, we explain why the Mac Lane cohomology groups
of Definition \ref{HML_gr_def} provide an appropriate tool for the systematic study of Toda brackets.
\begin{definition} 
Let $\cT$ be a triangulated category with a compact object $N$ such that $\Lambda = \cT(N,N)_*$ is 
$n$-sparse. $F_{\cT}(N,n)$ is defined to be the full subcategory of $\cT$ given by finite sums of copies of $N$
which are shifted by integral multiples of $n$. 
\end{definition}

The functor $\cT(N,-)_*$ induces an equivalence between $F_{\cT}(N,n)$ and the category $F(\Lambda,n)$.
This equivalence induces an isomorphism between the Mac Lane-cohomology group 
$\HML_{n\spl}^{*,-n}(\Lambda)$ and the normalized cohomology of $F_{\cT}(N,n)$ with coefficients in 
$\cT(-,-)_n$

Suppose we are given a sequence of composable maps $(\lambda'_1, \dots, \lambda'_{n+2})$ in $F(\Lambda,n)$ 
with $\lambda_{i+1} \lambda_{i} = 0$ for all $i$. We define the Toda bracket $\tbracket{\lambda'_1, \dots, \lambda'_{n+2}}$ of this sequence of maps
{\em in $\modu{\Lambda}$} to be the Toda bracket of the sequence  $(\lambda_1, \dots, \lambda_{n+2})$ in 
$F_{\cT}(N,n)$ associated to it under 
the equivalence $\cT(N,-)_*$. If $\cT$ is the derived category of a dga $A$, this defines Toda brackets
in the cohomology ring $H^*(A)$ via the Toda brackets in the derived category $\cD(A)$. One can check that
this recovers the usual notion of Massey products. 

\begin{remark} \label{comparing_indeterminacies}
In the situation above, the indeterminacy of $\tbracket{\lambda_1, \dots, \lambda_{n+2}}$
is $(\lambda_1)_{*} (\cT(X_{n+2}[n], X_1)) + (\lambda_{n+2}[n])^* (\cT(X_{n+1}[n],X_0))$ by Proposition 
\ref{indet+existence}. Now suppose we are given a normalized cocycle $c$ representing a cohomology class
$\gamma \in H^{n+2}(F_{\cT}(N,n), \cT(-,-)_n)$. Then $c(\lambda_1, \dots, \lambda_{n+2}) \in \cT(X_{n+2}[n],X_0)$. 
If we change $c$ by adding a coboundary $\delta(b)$, the evaluation on $(\lambda_1, \dots, \lambda_{n+2})$ changes by an element of 
$(\lambda_1)_{*} (\cT(X_{n+2}[n], X_1)) + (\lambda_{n+2}[n])^* (\cT(X_{n+1}[n],X_0)).$ 

Hence the {\em evaluation
of a cohomology class} has the same indeterminacy as the $(n+2)$-fold Toda bracket. Consequently, it makes sense to ask 
the evaluation of a cohomology class $\gamma \in \HML^{n+2,-n}_{n\spl}(\Lambda)$ on a complex 
of $n$-split $\Lambda$-modules $(\lambda'_{1}, \dots, \lambda'_{n+2})$ to {\em be} the Toda bracket 
$\tbracket{\lambda'_1, \dots, \lambda'_{n+2}}$ without having to mention indeterminacies. 
In other words, the indeterminacy of Toda brackets is built
into the cohomology of categories. For $n=3$, this observation was used for the study of (triple) universal 
Toda brackets in \cite{Baues_D_89}.
\end{remark}

\begin{theorem}  \label{toda_b-obstr}
Let $\cT$ be a triangulated category, and let $N$ be a compact object such that $\Lambda = \cT(N,N)_*$ 
is $n$-sparse. Let $M$ be a $\Lambda$-module admitting  a resolution 
\[\dots \xrightarrow{\lambda'_1} M_0 
\xrightarrow{\lambda'_0} M \to 0\]
by finitely generated free $n$-sparse $\Lambda$-modules. Let $\gamma \in \HML^{n+2,-n}_{n\spl}(\Lambda)$ 
be a cohomology class  
such that the evaluation $\gamma(\lambda'_1, \dots, \lambda'_{n+2})$
is the Toda bracket $\tbracket{\lambda'_1, \dots, \lambda'_{n+2}}$. Then the product 
$\id_M \cup \gamma  \in \Ext_{\Lambda}^{n+2,-n}(M,M)$
coincides with the unique obstruction class $\kappa_{n+2}(M)$ of Corollary \ref{uniqueness_obstr}.
\end{theorem}

\begin{proof}
We denote the realization of the resolution of $M$ by $N$-free objects by 
\[X_{n+2} \xrightarrow{\lambda_{n+2}} X_{n+1} \xrightarrow{\lambda_{n+1}} \dots 
 \xrightarrow{\lambda_1} X_0,\] 
so $(\lambda_i)_* = \lambda'_i.$
By Lemma \ref{exis_filtered}, there is a unique $n$-filtered object $Z \in \{\lambda_2, \dots, \lambda_{n}\}$. 
Since the $(n+1)$-fold Toda bracket of $(\lambda_1, \dots, \lambda_{n+1})$ contains only zero for degree reasons,
we can find maps $\alpha \colon Z \to X_0$ and $\beta \colon X_{n+1}[n-1] \to Z$ with $\sigma_Z \beta = \lambda_{n+1}[n-1]$ and $\alpha \sigma'_Z = \lambda_1$ such that $\alpha \beta = 0 \in \tbracket{\lambda_1, \dots, \lambda_{n-1}}$.

We use $\alpha$ and $\beta$ to find
distinguished triangles 
\[
Z \xrightarrow{\alpha} X_0 \to Y \xrightarrow{\omega} Z[1] \qquad \text{ and } \qquad  
X_{n+1}[n-1] \xrightarrow{\beta} Z \xrightarrow{\iota} X \to X_{n+1}[n].
\]
Lemma \ref{cones_filtered} tells us that $X$ is an $(n+1)$-filtered object in $\{\lambda_2, \dots, \lambda_{n+1}\}$ and
that $Y$ is an $(n+1)$-filtered object in $\{\lambda_1, \dots, \lambda_n\}$. 

The Toda bracket of $(\lambda_1, \dots, \lambda_{n+2})$ is non-empty by Proposition \ref{indet+existence}.
It can be defined using the $n$-filtered object $X$. 
Hence there are maps $\gamma_0 \colon X \to X_0$ and $\gamma_{n+2} \colon X_{n+2}[n] \to X$ with
$\gamma_0 \sigma'_X = \lambda_1$ and $\sigma_X \gamma_{n+2} = \lambda_{n+2}[n]$ 
such that $\gamma' = \gamma_0 \gamma_{n+2}$ is an element of $\tbracket{\lambda_1, \dots, \lambda_{n+2}}$. 
Looking at the triangle defining $X$, we see that $\gamma_0$ can be constructed by extending 
$\alpha \colon Z \to X_0$ to a map $X \to X_0$. The relation $\gamma_0 \iota = \alpha$ implies the existence of
the map $\rho$ in the following commutative diagram:
\[ \xymatrix{& & X_{n+2}[n] \ar[d]^{\lambda_{n+2}[n]} \ar[dl]_{\gamma_{n+2}} & \\
Z \ar[d]_{=} \ar[r]^{\iota} & X \ar[d]_{\gamma_0} \ar[r] & X_{n+1}[n] \ar[r]^{\beta[1]} \ar@{-->}[d]_{\rho} & Z[1] \ar[d]^{=}\\
Z \ar[r]^{\alpha} & X_0 \ar[r]^{\sigma'_Y} & Y \ar[r]^{\omega} & Z[1].}
\]
Here we use that the map $X_0 \to Y$ from the distinguished triangle defining $Y$ coincides with the map 
$\sigma'_Y$ which is part of the data of the $n$-filtered object $Y$. 

Applying $\cT(N,-)_*$ to the last diagram, we obtain the following commutative diagram of $\Lambda$-modules: 
\[ \xymatrix@-0.5pc{
\cT(N, X_{n+2}[n])_*  \ar[d]_{\gamma'_*} \ar[rr]^{(\lambda_{n+2})_*} & & \cT(N,X_{n+1}[n])_* \ar[d]_{\rho_*} 
\ar[rr]^{(\lambda_{n+1})_*} & &\cT(N, X_n[n])_* \ar@{=}[d] \ar[r] & \dots \\
\cT(N,X_0)_* \ar[rr]^{(\sigma'_Y)_*} \ar[dr]_{\lambda'_0} & & \cT(N,Y)_* \ar[rr]^{((\sigma_Z[1])\omega)_*} 
& & \cT(N,X_n[n])_* \ar[r] & \dots \\
& M \ar[ur] 
} \]
The lower sequence starting with $M$ in this diagram represents $\id_M \cup \gamma$ up to sign. 
Inspecting (\ref{representing_obstr_cls}) and  Lemma \ref{filtered_Postnikov}, we observe that 
it, up to signs, represents as well the exact sequence associated
to the $(n+1)$-Postnikov system obtained from $Y$. This uses that the map $(\sigma_Z[1])\omega$
equals the map $p_{n+1}$ of the $(n+1)$-filtered object $Y$, and therefore the map $(-1)^n \pi_n[n]$ of the associated 
Postnikov system. The sign of the latter map cancels with the $n$ factors  $(-1)$ by which the maps $(\lambda_i)_*$ 
differ from the differentials of the resolution induced by the Postnikov system. The remaining sign $(-1)^{\frac{n(n+3)}{2}}$
of the map $\sigma'_Y$ cancels with the sign built into the cup product.  
\end{proof}

Applications of this theorem will be given in Section \ref{ring_sp_section}. We point out that for $n=1$, the last theorem 
also leads to an interpretation of the product of a $\Lambda$-module homomorphism $f \colon M \to M'$ with $\gamma$,  
provided that $M$ satisfies the hypothesis of the theorem: by \cite[Proposition 3.4(iv) and Theorem 3.7]{Benson_K_S_04}
and the naturality of the cup product, $f \cup \gamma$ vanishes if and only if $f$ factors through a realizable
$\Lambda$-module. 

\section{Construction of universal Toda brackets} \label{constr_section}
As outlined in the introduction, the characteristic Hochschild cohomology class $\gamma_A$ 
of a dga $A$ considered in \cite{Benson_K_S_04}
is a motivation for the study of the universal Toda bracket $\gamma_R$ of a ring spectrum $R$. The class $\gamma_A$ cannot be recovered
from the derived category $\cD(A)$ \cite[Example 5.15]{Benson_K_S_04}. This suggests that the construction of $\gamma_R$ from 
$R$ will need more input than $\Ho(\modu{R})$. 
It turns out that the stable model structure on the category $\modu{R}$ together with the topological enrichment provides
the necessary information.  

Having the example $\modu{R}$ in mind, 
we construct the universal Toda bracket of an $n$-split subcategory of a general stable topological model category
in this section. The applications to ring spectra and the link to the realizability obstructions discussed above
are given in Section \ref{ring_sp_section}.

Besides \cite{Benson_K_S_04},  Baues' work on universal triple Toda
brackets \cite{Baues_97,Baues_D_89} is another motivation for our construction (and its name). 
He is working mainly in an unstable context,
considering  subcategories of $H$-group or $H$-cogroup objects in the homotopy category of 
topological spaces, though he 
points out that these constructions generalize to `cofibration categories' 
\cite[Remark on p. 271]{Baues_97}. We will only work in a stable context, in order to provide the
link to triangulated categories. This also avoids certain difficulties in the unstable case arising from maps which
are not suspensions (see the correction of \cite{Baues_D_89} in \cite[Remark on p. 270]{Baues_97}). 
We also do not use Baues' language of
`linear track extensions', as these seem to be only appropriate for the study of triple
universal Toda brackets. Nevertheless, the $n=1$ case of   
Proposition \ref{discard_basepoint} is basically what Baues encodes in a linear track extension. 

A motivation for the actual construction of the representing
cocycle is the approach of Blanc and Markl to higher homotopy operations \cite{Blanc_M_2003}. 
For a directed category $\Gamma$, the authors use the bar resolution $W \Gamma$ in the sense
of Boardman and Vogt \cite[III, \S 1]{Boardman_V_1973} to define general 
higher homotopy operations. If $\Gamma$ is the category generated by $n+2$
composable morphisms, this specializes to the higher Toda brackets we would
like to construct. In this case, $W \Gamma$ is just an $(n+1)$-dimensional 
cube. As we are not interested in other indexing categories, we will just use 
the cubes and do not make use of the bar resolution in our construction. 

In what follows, we assume familiarity with model categories. Hovey's book \cite{Hovey-model_1999}
provides a good reference. Other than in Quillen's original treatment of model categories \cite{Quillen-homotopical_67}, 
we will follow Hovey in assuming our model categories to have all small limits and colimits as well as functorial
factorizations. $\Top$ will be the category of compactly generated weak Hausdorff spaces, and $\Top_*$ will be the
pointed version.  The reason for working with these categories of spaces is
that $\Top_*$  is a closed symmetric monoidal model category \cite[Corollary 4.2.12]{Hovey-model_1999}.
We will often use stable topological model categories that are built on $\Top_*$. 
See Appendix \ref{top_model_cat} for a brief review.

\subsection{Cube systems}
Some notation is needed to state the next definition. Let $N(\cU)$ be the nerve of a small category $\cU$.
We write $d_i \colon N_{n}(\cU) \to N_{n-1}(\cU)$ for the $i$th simplicial face map, 
and $d_i^{\textrm{fr}} \colon N_{n}(\cU) \to N_i(\cU)$ and $d_{i}^{\textrm{ba}} \colon N_{n}(\cU) \to N_{i}(\cU)$ 
for the simplicial `front face' and the `back face' maps. 
In our notation for sequences of composable maps
from the preceding sections, this means for example $d_{n-i}(f_1, \dots, f_n) = (f_1, \dots, f_i f_{i+1}, \dots, f_{n})$ and  
$d_i^{\textrm{ba}}(f_1, \dots, f_n)=(f_1, \dots, f_i)$. 
We resist from reversing the notation for $(f_1, \dots, f_n)$ to make these formulas more intuitive here, since this would
be inconsistent with our previous convention, which was chosen since $(f_1, \dots, f_n)$ frequently arose from a projective 
resolution.  

For $n \geq 1$, we denote the $n$-fold cartesian product of 
the unit interval by $\cub{n}$ and define  $\cub{0}$ to be the one point space. 
For $\epsilon \in \{0,1\}$ and $1 \leq i \leq n$,
we have a structure map 
\[\epsilon^i \colon \cub{n-1} \to \cub{n}, \quad (t_1, \dots, t_{n-1}) \mapsto (t_1, \dots, t_{i-1}, \epsilon, t_{i}, \dots, t_{n-1}).\]
With $\epsilon, \omega \in \{0, 1\}$, these maps satisfy the relation $\epsilon^i \omega^{j-1} = \omega^j \epsilon^{i}$ if $1 \leq i < j \leq n$.
We write $\skel_{i}\cub{n}$ for the $i$-skeleton of $\cub{n}$ in the obvious CW-structure. 
When we consider $\cub{n}$ with $n \geq 1$ as a pointed space, we take $(1,
\dots, 1)$ as the basepoint.

For a stable topological model category $\cC$, we will work with the set of maps 
$\Top(\cub{n}, \Map_{\cC}(X,Y))$. The enriched composition $\mu$ of $\cC$ induces a composition
\[\begin{split}
\mu_{p,q} \colon & \Top(\cub{p}, \Map_{\cC}(Y,Z)) \times \Top(\cub{q}, \Map_{\cC}(X,Y)) \to \Top(\cub{q+p}, \Map_{\cC}(X,Z)), \\ 
&(b, b') \mapsto ((t_1, \dots, t_{p+q}) \mapsto (x \mapsto  b_{(t_1, \dots, t_p)}(b'_{(t_{p+1}, \dots, t_{p+q})}(x)))).
\end{split}\]
The associativity of the enriched composition implies that $\mu_{p,q+r}(\id \times \mu_{q,r})$ 
and $\mu_{p+q, r}(\mu_{p,q} \times \id)$ correspond under the coherence isomorphism for associativity
of the $3$-fold cartesian product in $\Top$. 

The zero map is a canonical basepoint for $\Map_{\cC}(X,Y)$. When a possibly different
map $g \colon X \to Y$ is used as the basepoint, we write $(\Map_{\cC}(X,Y),g)$ for the resulting pointed space. 

For $\epsilon \in \{0,1\}$ and $1 \leq i \leq p$, we have $(\epsilon^i)^* \colon \Top(\cub{p}, T) \to \Top(\cub{p-1},T).$ 
If $1 \leq i \leq p+q$, these restrictions satisfy
\[ (\epsilon^i)^* \mu_{p,q}(b,b') = \begin{cases} \mu_{p-1,q}((\epsilon^i)^*b, b') & \quad \textrm{if } i \leq p, \\
\mu_{p,q-1}(b,(\epsilon^{i-p})^*b') & \quad \textrm{if } i > p. \end{cases}\]

\begin{definition}  \label{cube_system_def}
Let $\cU$ be a small full subcategory of the homotopy category of a stable topological model category $\cC$. A cube 
system for $\cU$ consists of the following data: for every object $X$ of $\cU$, there is a 
cofibrant and fibrant object $\Phi(X)$ of $\cC$ and an isomorphism $\varphi_X \colon X \to \Phi(X)$
in $\Ho(\cC)$. We write $\Phi(\cU)$ for the set of all those objects. Furthermore, for $0 \leq j \leq n$
there are maps 
\[ b^j \colon N_{j+1}(\cU) \to \coprod_{X,Y \in \Phi(\cU)} \Top(\cub{j}, \Map_{\cC}(X,Y)) \]
such that 
\begin{enumerate}[(i)]
\item $b^0(f \colon X \to Y) \in \Map_{\cC}(\Phi(X), \Phi(Y))$ and $b^0(f)$ represents $f$ on the
model category level, i.e., $\varphi_Y^{-1} b^0(f) \varphi_X = f$ in $\Ho(\cC)$. 
\item If one of the maps $f_i$ in $(f_1, \dots, f_{j+1})$ is a zero map, then $b^j(f_1, \dots, f_{j+1})$
has the zero map in $\cC$ as constant value. 
\item 
For $j \geq 1$ and $1 \leq i \leq j$, $b^{j-1} d_{j+1-i} = (1^i)^* b^j$. 
\item 
For $j \geq 1$ and $1 \leq i \leq j$, 
\[\mu_{i-1,j-i}(b^{i-1} \times b^{j-i}) (d_i^{\textrm{ba}} \times d_{j+1-i}^{\textrm{fr}}) \Delta = (0^i)^* b^j,\]
where $\Delta$ is the diagonal and $\mu_{i-1,j-i}$ is explained above. 
\end{enumerate}
By (iii), $b^j(f_1, \dots, f_{j+1})$ maps the basepoint of $\cub{j}$
to $b^0(f_1 \dotsm f_{j+1})$. 
\end{definition}

\begin{figure}[b]
\[\xymatrix@-.25pc@!R{
& & & (f_1)(f_2 f_3)(f_4) \ar@{-}[dddlll]_(.35){(f_1)\circ(f_2f_3)(f_4)} \ar@{-}[rrr]^(.52){(f_1)(f_2)(f_3)\circ(f_4)} \ar@{--}[dd] 
& & &  (f_1)(f_2)(f_3)(f_4) 
\ar@{-}[dd] \ar@{-}[dddlll]_(.35){(f_1)\circ(f_2)(f_3)(f_4)} \\
& & & & & & \\
& & & (f_1)(f_2f_3f_4) \ar@{--}[dddlll]^(.65){(f_1)\circ(f_2f_3f_4)} \ar@{--}[rrr]_(.52){(f_1)(f_2)\circ(f_3f_4)} & & & 
(f_1)(f_2)(f_3f_4) \ar@{-}[dddlll]^(.65){(f_1\circ f_2)(f_3f_4)} \ar@{-} \\
{(f_1f_2f_3)(f_4)}  \ar@{-}[rrr]^(.45){(f_1 f_2)\circ( f_3)(f_4)}  & & & (f_1f_2)(f_3)(f_4)  \ar@{-}[dd] & & &  \\
 & & & & & & \\
(f_1 f_2 f_3 f_4) \ar@{-}[rrr]_(.45){(f_1f_2)\circ(f_3f_4)} \ar@{-}[uu]   & & & (f_1 f_2)(f_3 f_4) } \]
\caption{A $3$-cube.} \label{three_cube}
\end{figure}
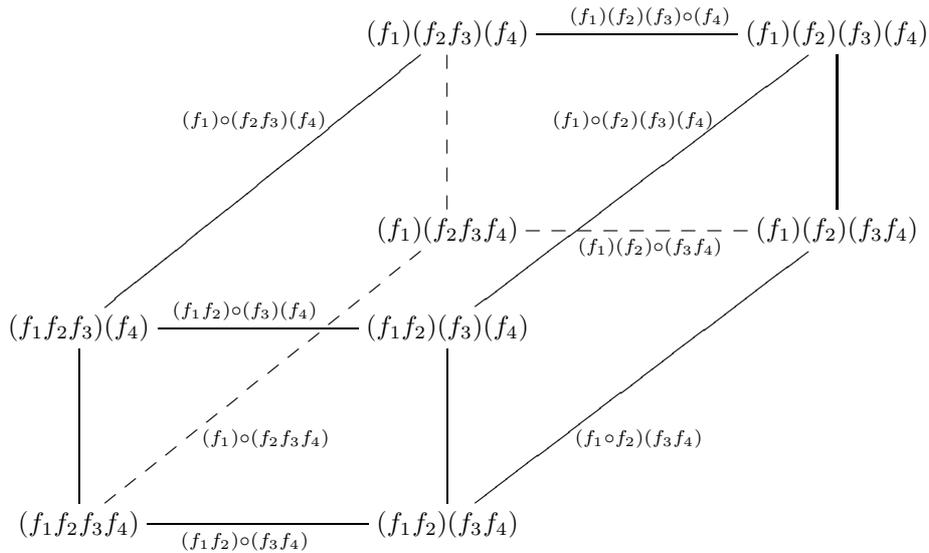

A $0$-cube system chooses maps in the model category representing maps in the homotopy category. In general, 
it is not possible to arrange these choices such that $b^0(f_1)b^0(f_2) = b^0(f_1 f_2)$ holds. Nevertheless, 
these maps are homotopic, and unraveling (ii) and (iii) shows  that a $1$-cube system specifies a homotopy 
 $b^1(f_1,f_2)$ between them. For $j \geq 2$, the $b^j(f_1, \dots, f_{j+1})$ encode coherence homotopies between
different choices of representatives and coherence homotopies of lower degree. 
Figure \ref{three_cube} 
(compare \cite[Figure 2.12]{Blanc_M_2003})
illustrates the case $n=3$. In the picture, we write $(f_j \dotsm f_k)$ for $b^0(f_j \dotsm f_k)$ and 
$(f_j \dotsm f_{k-1})\circ(f_k \dotsm f_{l})$ for $b^1(f_j \dotsm f_{k-1},f_{k} \dotsm f_{l})$. 

\begin{definition}  
In the situation of Definition \ref {cube_system_def}, a {\em pre $n$-cube system} for $\cU$ consists of 
an $(n-1)$-cube system for $\cU$ and a map 
\[ \widehat{b}^n \colon N_{n+1}(\cU) \to \coprod_{X,Y \in \Phi(\cU)} \Top(\skel_{n-1} \cub{n}, \Map_{\cC}(X,Y)) \]
such that $\widehat{b}^n$ and the $b^j$ for $j<n$ satisfy conditions (ii)-(iv) of Definition \ref {cube_system_def}. This makes
sense since (iii) and (iv) only involve the behavior on $\skel_{n-1} \cub{n}$. 
\end{definition}
Similar as above, $\widehat{b}^n(f_1, \dots, f_{n+1})$ maps the basepoint to 
$b^0(f_1 \dotsm f_{n+1})$.

\begin{lemma} \label{extending_cube_systems}
An $(n-1)$-cube system for $\cU$ can be extended to a pre $n$-cube system. 
The restriction of $\widehat{b}^n(f_1, \dots ,f_{n+1})$ to the subcubes $(0^i)(\cub{n-1})$ for $1 < i < n$ is 
determined by the underlying $(n-2)$-cube system. 
\end{lemma}
\begin{proof}
Since $\skel_{n-1}\cub{n}$ is the union of the $(n-1)$-dimensional subcubes $(0^i)(\cub{n-1})$ and 
$(1^i)(\cub{n-1})$, we define the restriction of $\widehat{b}^n$ to these subcubes by
\[ (1^i)^* \widehat{b}^n := b^{n-1} d_{n+1-i} \quad \textrm{ and } \quad (0^i)^* \widehat{b}^n := \mu_{i-1,n-i}(b^{i-1} \times b^{n-i}) (d_i^{\textrm{ba}} \times d_{n+1-i}^{\textrm{fr}}) \Delta.\]
 It remains to check that this is well defined on the intersections. 

Let $1 \leq j < k \leq n$. On 
$1^k(1^j(\cub{n-2})) = 1^j(1^{k-1}(\cub{n-2}))$ we
have 
\[(1^{k-1})^* (1^{j})^* \widehat{b}^n = b^{n-1} d_{n-(k-1)} d_{n+1-j} = b^{n-1} d_{n-j} d_{n+1-k} = (1^j)^* (1^{k})^* \widehat{b}^n.\] 
Next we check the compatibility on $1^k(0^j(\cub{n-2})) = 0^j(1^{k-1}(\cub{n-2})).$ 
A somewhat lengthy calculation involving the interchange formula for $\mu_{p,q}$ and
$(1^i)^*$ mentioned above shows that both $(1^{k-1})^* (0^j)^* \widehat{b}^n$
and $(0^j)^* (1^{k})^* \widehat{b}^n$ equal 
\[\mu_{j-1,n-j-1}(b^{j-1} \times b^{n-j-1})(d_j^{\textrm{ba}} \times d_{n-k+1} d_{n+1-j}^{\textrm{fr}}) \Delta.\] 

The case of $(0^{k-1})^* (1^j)^* \widehat{b}^n$
and $(1^j)^* (0^{k})^* \widehat{b}^n$ is similar. The remaining case of  $(0^{k-1})^* (0^j)^* \widehat{b}^n$
and $(0^j)^* (0^{k})^* \widehat{b}^n$ follows from the associativity of the maps $\mu_{p,q}$. 
\end{proof}

\begin{proposition} \label{cube_system_existence}
Let $\cU $ be a small $n$-split subcategory of the homotopy category of a stable topological
model category $\cC$. Then there exists an $n$-cube system for $\cU$. 
\end{proposition}
\begin{proof}
We start with choosing the representing objects $\Phi(X)$, the isomorphisms $\varphi_X$ and the 
representing maps $b^0(f)$, which is always possible. Here we consider $b^0(f)$ as an element
of $\Top(\cub{0},\Map_{\cC}(\Phi(X), \Phi(Y)))$. For each pair of composable maps $(f_1,f_2)$, 
the maps $b^0(f_1 f_2)$ and $b^0(f_1)b^0(f_2)$ are homotopic. After adjoining and forgetting 
the basepoint, a homotopy $(\cub{1})_+ \sm \Phi(X_2) \to \Phi(X_0)$ gives rise to $b^1(f_1, f_2)$. 
This completes the $1$-cube system. 

Inductively, suppose we have constructed a $j$-cube system for $j<n$. By Lemma \ref{extending_cube_systems},
it induces a pre $(j+1)$-cube system. In order to extend $\widehat{b}^{j+1}(f_1, \dots, f_{j+2})$ from 
$\skel_j \cub{j+1}$ to $\cub{j+1}$, it suffices to know that it represents the trivial homotopy class in 
$\pi_j(\Map_{\cC}(\Phi(X_{j+2}),\Phi(X_{0})), b^0(f_1 \dotsm f_{j+2}))$. This group is isomorphic to 
$[S^j \sm X_{j+2}, X_0]^{\Ho(\cC)} \iso [X_{j+2},X_0]^{\Ho(\cC)}_j$ by means of the isomorphism $\sigma_{(f_1 \dotsm f_{j+2})}$
of Proposition  \ref{discard_basepoint}, and hence trivial since $\cU$ is $n$-split. 
\end{proof}

\subsection{The universal Toda bracket} We now put the data of a cube system together
to get the desired cohomology class. 
\begin{construction} \label{gen_construction}
Let $\cC$ be a stable topological model category and let $\cU$ be a small $n$-split subcategory
of $\Ho(\cC)$. Then there is a well defined cohomology class  
$\gamma_{\cU} \in H^{n+2}(\cU, [-,-]^{\Ho(\cC)}_n)$ 
which determines by evaluation all $(n+2)$-fold Toda brackets
of complexes of $n+2$ composable maps in $\cU$. 

We choose an $n$-cube system for $\cU$ which is possible by Proposition \ref{cube_system_existence},
and extend it to a pre $(n+1)$-cube system by Lemma \ref{extending_cube_systems}. Then we define 
a normalized cochain $c \in \ovl{C}^{n+2}(\cU, [-,-]_n^{\Ho(\cC)})$ as follows. Its evaluation 
on a sequence of \mbox{$(n+2)$} composable maps 
$ X_{n+2} \xrightarrow{f_{n+2}} X_{n+1} \xrightarrow{f_{n+1}} \dots \xrightarrow{f_1} X_0 $
in $\cU$ 
is the image of the homotopy class of the (pointed) map $\widehat{b}^{n+1}(f_1, \dots, f_{n+2})$ under the chain of isomorphisms
\[
\xymatrix@-.25pc{
[\skel_{n} \cub{n+2}, (\Map_{\cC}(\Phi(X_{n+2}), \Phi(X_{0})), b^0(f_1 \dotsm f_{n+2}))]^{\Ho(\Top_*)}  
\ar[d]_{\sigma_{(f_1 \dotsm f_{n+2})}} \\ 
\qquad \qquad \qquad \qquad [\skel_{n}\cub{n+1} \sm X_{n+2}, X_0 ]^{\Ho(\cC)} \xrightarrow{\iso} [X_{n+2}[n],X_0]^{\Ho(\cC)}} \]

We show in Lemma \ref{construction_yields_cocycle} that $c$ is a cocycle. Lemma \ref{construction_well_defined} verifies that
its cohomology class does not depend on the choice of the cube system. In Proposition \ref{Tbrackets_comp}, we show
that the evaluation of $c$ on a complex of maps in $\cU$ is an element of the Toda bracket of that complex. With the 
comparison of the indeterminacies in Remark \ref{comparing_indeterminacies}, it follows that 
the evaluation of the {\em cohomology class} 
$\gamma_{\cU}$ of $c$ on a complex in $\cU$ yields its Toda bracket. This is why we call the well defined cohomology 
class $\gamma_{\cU}$  the {\em universal Toda bracket} of $\cU$. 
\end{construction}

The Homotopy Addition Theorem \cite{Hu_1953} will be our tool in the proof of the next two lemmas. We use
$T=\{(\epsilon, i) | \epsilon \in \{0,1\}, 1 \leq i \leq n+2 \}$ as indexing set for the $(n+1)$-dimensional
subcubes of $\cub{n+2}$. It is a disjoint union of $T_{+} = \{(\epsilon, i) | (-1)^{n+\epsilon + 1 } = 1 \}$
and $T_{-} = \{(\epsilon, i) | (-1)^{n+\epsilon + 1 } = -1 \}$. 

Given a pointed topological space $K$ with abelian fundamental group and a pointed map $f \colon \skel_n \cub{n+2} \to K$,
the Homotopy Addition Theorem states that
\[ \sum_{(\epsilon, i) \in T_{+}} [ (\epsilon^i)^* f ] - \sum_{(\epsilon, i) \in T_{-}} [ (\epsilon^i)^* f ] = 0. \]
Here $[(\epsilon^i)^* f ]$ is the homotopy class of the restriction of $f$ to $(\epsilon^i)^*(\skel_n \cub{n+2})$, 
or the image of the operation of a path to the basepoint of $\skel_n \cub{n+2}$ on $[(\epsilon)^* f]$ if 
$(\epsilon^i)^*(\skel_n \cub{n+2})$ doesn't contain the basepoint of $\skel_{n} \cub{n+2}$. This is well defined since
two such paths are homotopic if $n>1$ and the $\pi_1$-action is trivial for $n=1$ as $\pi_1$ is assumed to be abelian. 
Our source for this formulation of the theorem is \cite[VII.9.6]{Bredon-topology_1997}. The signs result from 
specifying an orientation through choosing `even' and `odd' faces $T_{+}$ and $T_{-}$, compare 
\cite{Whitehead-homotopy_1953}.

\begin{lemma} \label{construction_yields_cocycle}
The cochain $c$ of Construction \ref{gen_construction} is a cocycle.
\end{lemma}

\begin{proof}
We fix a sequence of $(n+3)$ composable maps $(f_1, \dots, f_{n+3})$ in $\cU$.  
As in Lemma \ref{extending_cube_systems}, the pre $(n+1)$-cube system induces 
\[e(f) = e(f_1, \dots, f_{n+3}) \colon \skel_{n} \cub{n+2} \to \Map_{\cC}(\Phi(X_3), \Phi(X_0))\] with
$(1^i)^* e(f):= (\widehat{b}^{n+1} d_{n+3-i})(f_1, \dots, f_{n+3})$ and 
\[ (0^i)^* e(f) := (\mu_{i-1,n+2-i}(\widehat{b}^{i-1} \times \widehat{b}^{n+2-i}) (d_i^{\textrm{ba}} \times d_{n+3-i}^{\textrm{fr}}) \Delta)(f_1, \dots, f_{n+3}), \]
where $\widehat{b^i}$ is the restriction of $b^i$ if $i<{n+1}$. 

We apply the Homotopy Addition Theorem mentioned above to this map to get 
\[ 0= \sum_{(\epsilon, i) \in T_{+}} [ (\epsilon^i)^* e(f)] - 
\sum_{(\epsilon, i) \in T_{-}}  [ (\epsilon^i)^* e(f)] \]
in $\pi_n(\Map_{\cC}(\Phi(X_{n+3}), \Phi(X_0)), b^0(f_1 \dotsm f_{n+3}))$. For $n=1$, this uses that
$\cC$ being stable implies $\pi_1$ is abelian. For $1 < i < n+2$, the restrictions $(0^i)^* e(f)$
extend to maps $(0^i)^* (\skel_{n+1} \cub{n+1})$ as $\widehat{b}^{i-1} = b^{i-1}$ and $\widehat{b}^{n+2-i} = b^{n+2-i}$. Hence their homotopy classes vanish. 
Now we apply the isomorphism $\sigma_{(f_1 \dotsm f_{n+3})}$ to the sum and use that $\sigma$ is additive and invariant 
under the action of basepoint changing paths. Hence 
\[ \begin{split}
0 & = (f_1)_* c(f_2, \dots, f_{n+3}) + \sum_{i=1}^{n+2} (-1)^i c(f_1,\dots, f_i f_{i+1}, \dots, f_{n+3}) \\
& + (-1)^{n+3} (f_{n+3})^* c(f_1, \dots, f_{n+2}).
\end{split} \]
\end{proof}

\begin{lemma} \label{construction_well_defined}
The cohomology class $\gamma_{\cU}$ of Construction \ref{gen_construction} does 
not depend on the choice of a cube system.  
\end{lemma}
\begin{proof}
In the first step we assume we are given a second $(n-1)$-cube system $(\udl{\Phi}, \udl{\varphi}, \udl{b^j})$ for $\cU$. 
We show that it extends to an $n$-cube system giving the same cohomology class. 

For every object $X$ in $\cU$, 
our data specifies an isomorphism  $\Phi(X) \to \udl{\Phi}(X)$ in $\cU$. We realize it by $g_X$ in $\cC$ and its inverse
by $g'_X$. 
For $f_1 \colon X_1 \to X_0$ in $\cU$, we know $g_{X_0} b^0(f_1) g'_{X_1} = \udl{b}^0(f_1)$ in $\cU$. Let $h^0(f) \colon \cub{1} \to
\Map_{\cC}(\udl{\Phi}(X_1), \udl{\Phi}(X_0))$ be a homotopy between them. 
With similar arguments as in Proposition 
\ref{cube_system_existence}, one can iterate the construction to find maps 
\[h^j \colon N_{j+1}(\cU) \to \coprod_{X, Y \in \cU} \Top(\cub{j+1}, \Map_{\cC}(\udl{\Phi}(X), \udl{\Phi}(Y)))\]
for $j < n$ with $(0^1)^* h^j = (g_{X_0})_* (g'_{X_{j+1}})^* b^j$ and $(1^1)^* h^j =  \udl{b}^j$.
For $j=n$, we use the homotopy extension property of $(\skel_n \cub{n+1}) \setminus (1^1)(\cub{n}) \to \cub{n+1}$
to find an $h^n$ whose restriction to $(1^1)(\cub{n})$ defines $\udl{b}^n$. When we form the pre $(n+1)$-cube systems associated
to the two cube systems, the $h^j$ assemble to a homotopy between $(g_{X_0})_* (g'_{X_{j+2}})^* \widehat{b}^{n+1}$ and $\udlhat{b}{n+1}$. Hence $\sigma$ associates the same cocycle to them.

Now suppose we are given two $n$-cube systems  for $\cU$. The first part shows that we can assume their underlying $(n-1)$-cube systems
to coincide. Let $T_{\textrm{dev}} = \{(1,i)|1 \leq i \leq n+1\} \cup \{(0,1), (0,n+1)\}$. Lemma \ref{extending_cube_systems} shows that
the associated pre $(n+1)$-cube systems can only deviate on the faces specified by $T_{\textrm{dev}}$.

Let $A$ be the space obtained from $\skel_n{\cub{n+1}}$ by gluing for each $(\epsilon, i) \in T_{\textrm{dev}}$ one copy of $\cub{n} \cup_{\skel_{n-1} \cub{n}} \cub{n}$ along the right hand side copy of $\cub{n}$ to $(\epsilon^i)(\cub{n})$. Let $i \colon \skel_n \cub{n+1} \to A$ 
be the canonical inclusion and let $\udl{i} \colon \skel_n \cub{n+1} \to A$ be the injection which maps $(\epsilon^i)\cub{n}$ to 
the left copy of $\cub{n}$ in the pushout. Then the two pre cube systems induce a map 

\[\widehat{a} \colon N_{j+2}(\cU) \to \coprod_{X, Y \in \cU} \Top(A, \Map_{\cC}(\udl{\Phi}(X), \udl{\Phi}(Y)))\]
with $i^* \widehat{a} = \widehat{b}^{n+1}$ and $\udl{i}^* \widehat{a} = \udlhat{b}{n+1}$. By the slightly different incarnation 
of the Homotopy Addition Theorem  \cite[VII.9.5]{Bredon-topology_1997},
the evaluations on $(f_1, \dots, f_{n+2})$ satisfy

\[ \udlhat{b}{n+1} = 
\widehat{b}^{n+1} + \sum_{(\epsilon,i) \in T_{\textrm{dev}} \cap T_{+}} [(\iota_{(\epsilon,i)})^* \widehat{a}] - \sum_{(\epsilon,i) \in T_{\textrm{dev}} \cap T_{-}}[(\iota_{(\epsilon,i)})^* \widehat{a}], \]
where $\iota_{(\epsilon,i)} \colon  \cub{n} \cup_{\skel_{n-1} \cub{n}} \cub{n} \to A$ is the inclusion which belongs to 
$(\epsilon, i) \in T_{\textrm{dev}}$. The signs arise in the same way as in the last lemma.

Let $a \colon N_{n+1}(\cU) \to \coprod_{X, Y \in \cU} \Top( \cub{n} \cup_{\skel_{n-1} \cub{n}} \cub{n}, \Map_{\cC}(\udl{\Phi}(X), \udl{\Phi}(Y)))$ be the map which is $b^n$ one the right and $\udl{b}^n$ on the left copy of $\cub{n}$. Applying $\sigma_{b^0(f_1 \dotsm f_{n+2})}$ to $[a(f_1, \dots, f_{n+1})]$ defines an $(n+1)$-cochain $\ovl{a} \in  \ovl{C}^{n+1}(\cU, [-,-]_n^{\Ho(\cC)})$. Next we apply $\sigma_{b^0(f_1 \dotsm f_{n+2})}$ to the sum formula above to get the desired equation 
\[ c(f_1, \dots, f_{n+2}) = \udl{c}(f_1, \dots, f_{n+2}) + (\delta \ovl{a})(f_1, \dots, f_{n+2}). \]
As in Lemma \ref{construction_yields_cocycle}, the orientations of the subcubes imply the signs needed for the coboundary formula.  
\end{proof}

The last lemma completes the proof of $\gamma_{\cU}$ being well defined. For later use, we prove two more lemmas closely
related to this construction. 

\begin{lemma} \label{morita_cubesys}
Let $G \colon \cC \to \cD$ be a left Quillen functor between stable topological model categories $\cC$ and $\cD$ which is compatible
with the topological structure. 
If $\cU$ and $\cW$ are 
small $n$-split subcategories of $\Ho(\cC)$ and $\Ho(\cD)$ such that $G$ induces an equivalence $\cU \to \cW$
and an isomorphism $G^*([-,-]_n^{\Ho(\cD)}) \iso [-,-]_n^{\Ho(\cC)}$, 
then the induced isomorphism
\[ G^* \colon H^{n+2}(\cW,[-,-]_n^{\Ho(\cD)}) \to H^{n+2}(\cU, [-,-]_n^{\Ho(\cC)}) \]
sends $\gamma_{\cW}$ to $\gamma_{\cU}$. 
\end{lemma}
\begin{proof}
We apply $G$ to the data of an $n$-cube system for $\cU$. This gives almost an $n$-cube system for $\cW$. The
only missing part is that the objects $G(\Phi(X))$ are not necessarily fibrant. 
Similarly as in Lemma \ref{construction_well_defined}, one can construct a cube system for $\cW$ such that
the resulting cocycle representing $\gamma_{\cW}$ becomes,  after
applying $G^*$, equivalent to that of the cube system for $\cU$. 
\end{proof}

\begin{lemma} \label{univ_tbracket_zero}
Let $\cC$ be a stable topological model category and let $\cU$ be a small
$n$-split subcategory of $\Ho(\cC)$. If $\gamma_{\cU}$ is trivial, then the map $b^n$ of any $n$-cube system for $\cU$ can
be changed such that the resulting new $n$-cube system has the zero cochain
as a representing cocycle. 
In particular, the modified $n$-cube system can be extended to an $(n+1)$-cube system. 
\end{lemma}
\begin{proof}
Since $\gamma_{\cU} = 0$, there is an $e \in \ovl{C}^{n+1}(\cU,[-,-]_n^{\Ho(\cC)})$ with $\delta (e) = c$. We model 
the $n$-sphere by gluing two copies of the $n$-cube $\cub{n}$ together along
their boundaries. 
For every sequence of composable maps $(f_1, \dots, f_{n+1})$ in $\cU$, there is a map 
$\widehat{e}(f_1, \dots, f_{n+1}) \colon (\cub{n} \cup_{\partial \cub{n}} \cub{n}) \to \Map_{\cC}(X_{n+1},X_0)$ such
that its restriction to the left copy of $\cub{n}$ is $b^n(f_1, \dots, f_{n+1})$ and that 
$[\sigma_{b^0(f_1 \dotsm f_{n+1})}(\widehat{e}(f_1, \dots, f_{n+1}))] = e(f_1, \dots , f_{n+1})$. We define 
$\udl{b}^n(f_1, \dots, f_{n+1})$ to be the restriction of $\widehat{e}(f_1, \dots, f_{n+1})$ to the right copy of $\cub{n}$. 
Similarly as in Lemma \ref{construction_well_defined}, the Homotopy Addition Theorem shows the assertion. 
\end{proof}

\section{Comparing definitions of Toda brackets}

This section is devoted to the proof of 
\begin{proposition} \label{Tbrackets_comp}
Let $\cC$ be a stable topological  model category, let $\cU$ be an $n$-split subcategory
of $\Ho(\cC)$, and let 
$X_{n+2} \xrightarrow{f_{n+2}} X_{n+1} \xrightarrow{f_{n+1}} \dots \xrightarrow{f_1} X_0 $
be a sequence of maps in $\cU$  with $f_{i}f_{i+1} = 0$ for $1 \leq i \leq n+1$.
Let $c$ be the cocycle of Construction \ref{gen_construction}.
Its evaluation on $(f_1, \dots, f_{n+2})$ lies in the Toda bracket $\tbracket{f_1, \dots, f_{n+2}}$ in the sense
of Definition \ref{higher_tbracket}.
\end{proposition}

\begin{remark} \label{discuss_different_tbrackets}
Triple Toda brackets were introduced by Toda \cite{Toda_1952,Toda-composition_1962} to study the stable homotopy
groups of spheres. Higher Toda brackets were introduced in the 1960's, and there are different approaches in the literature.
One of them is Cohen's definition using {\em filtered objects}
\cite[\S 2]{Cohen_1968}. We used a variant of this for triangulated categories in Definition \ref{higher_tbracket}.

Another approach is Spanier's definition of higher Toda brackets \cite{Spanier_1963} using the concept of a {\em carrier}. 
A related concept is Klaus' definition of a pyramid \cite[3.4]{Klaus_2001}, which is linked to Spanier's definition 
by \cite[Proposition 3.6]{Klaus_2001}. 
The perhaps most general approach to Toda brackets and other higher homotopy operations is that of Blanc 
and Markl \cite{Blanc_M_2003}, who define them as obstructions to realizing  homotopy commutative diagrams by 
strictly commutative ones. Their definition of Toda brackets  is related to Spanier's  \cite[Example 3.12]{Blanc_M_2003}. 

In Lemma \ref{evalua_descr} below we will see that the evaluation of the universal Toda bracket can be interpreted as something
similar to a pyramid in the sense of Klaus. Proposition \ref{Tbrackets_comp} shows that this is equivalent
to the Toda bracket defined via filtered objects in  Definition \ref{higher_tbracket}.
We work out the comparison as far as needed for our purposes in some detail since we were not able to find an appropriate reference 
in the literature which relates the different approaches. 
\end{remark}

For the rest of the section, we fix a $\cU$ and $(f_1, \dots, f_{n+2})$ with $f_i f_{i+1} = 0$ for $1 \leq i \leq n+1$ as in the proposition. We also fix an $n$-cube system $b^j$ of $\cU$ and write $\widehat{b}^{n+1}$ for the associated 
pre $(n+1)$-cube system and $c$ for the cocycle defined in Construction \ref{gen_construction}. 
For simplicity, we denote the objects $\Phi(X)$
of $\cC$ chosen by the cube system also by $X$. The map $\widetilde{b}^j(f_1, \dots, f_{j+1}) \colon (\cub{j})_+ \sm X_{j+1} \to X_0$
will always be the adjoint of the map $b^j(f_1, \dots, f_{j+1})$. 

We denote by  $\randcub{n+1}$ the pointed space obtained from $\skel_{n} \cub{n+1}$ by collapsing all $n$-dimensional 
subcubes $(1^i)(\cub{n})$ with $1 \leq i \leq n+1$ to the basepoint $(1, \dots, 1)$. 
The space $\randcub{n+1}$ is homeomorphic to an $n$-sphere. Since $f_i f_{i+1} = 0$, the map 
$\widehat{b}^{n+1}(f_1, \dots, f_{n+2})$
factors through the quotient map $\skel_n \cub{n+1} \to \randcub{n+1}$
and induces $\widetilde{c}^{n+1}(f_1, \dots, f_{n+2}) \colon \randcub{n+1} \to (\Map_{\cC}(X_{n+2},X_{0}),0)$.

\begin{lemma} \label{evalua_descr}
The map $\widetilde{c}^{n+1}(f_1, \dots, f_{n+2})$ represents the evaluation of $c$
on the complex $(f_1, \dots, f_{n+2})$ in $\cU$. 
\end{lemma}
\begin{proof}
This follows from Proposition  \ref{discard_basepoint}. 
\end{proof}

Depending on our chosen $n$-cube system and $(f_1, \dots, f_{n+2})$, we now construct
objects $F_j = F_j(f_2, \dots, f_{j+1})$ in $\cC$ for $ j \leq n+1$. Set 
\[A_j = A'_j = \coprod \limits_{1 \leq r < s \leq j+1} (\cub{j-1})_+ \sm X_{s} \qquad 
\text{and} \qquad B_j = \coprod \limits_{1 \leq i \leq j+1} (\cub{j})_+ \sm X_i. \]
The object $F_j$ is the coequalizer of two maps $h, k \colon A_j \coprod A'_j \to B_j$ to be described next.
For this, we think of the copies of $\cub{j}$ in $B_j$ as the $j+1$ subcubes $(0^i)(\cub{j})$ of $\cub{j+1}$. 
The copies of $\cub{j-1}$ in $A_j$ are thought of as the $(j-1)$-dimensional subcubes $(0^s)(0^r)(\cub{j-1})$ of $\cub{j+1}$, and 
the copies of $\cub{j-1}$ in $A'_j$ are thought of 
as the $(j-1)$-dimensional subcubes $(0^s)(1^r)\cub{j-1}$. 

The map $h$ is given by $(0^r)_+ \sm X_{s}$ on the copy of $(\cub{j-1})_+ \sm X_{s}$ in $A_j$ indexed by $(r,s)$, 
and by  $(1^r)_+ \sm X_{s}$ in the case of $A'_j$.  The map $k$ is the trivial map to the basepoint on $A'_j$. 
On the copy of $(\cub{j-1})_+ \sm X_{s}$ in $A_j$
indexed by $(r,s)$, it is given by the product of $0^{s-1}$ and the map $\widetilde{b}^{s-r-1}(f_{r+1}, \dots, f_s)$ using
the last $(s-r-1)$-coordinates of $\cub{j-1}$.

\begin{example}
The case $j=2$, which becomes relevant for $4$-fold Toda-brackets, is displayed in Figure \ref{fil_obj_fig}.
\begin{figure}
\[\xymatrix@-.2pc@!R{
& & & (0,1,0) \ar@{~}[ddddlll]|(.7){(\cub{1})_+ \sm X_3} \ar@{-}[rr]^{(\cub{1})_+ \sm X_3} 
\ar@{--}[dd]|{(\cub{1})_+ \sm X_1} & & (0,0,0) 
\ar@{-}[dd]|{(\cub{1})_+ \sm X_2} \ar@{-}[ddddlll]|(.7){(\cub{1})_+ \sm X_3} \\
&{}\save[]*{(\cub{2})_+ \sm X_1}\ar@<1ex>@/^1.5pc/@{.>}[rrr] \restore& & & & \\
& & & (0,1,1) \ar@{--}[rr]^{(\cub{1})_+ \sm X_1} & & (0,0,1) \ar@{-}[ddddlll]|(.7){(\cub{1})_+ \sm X_2} \\
& & & & & \\
(1,1,0) \ar@{~}[rr]_{(\cub{1})_+ \sm X_3}  & & (1,0,0) \ar@{~}[dd]|{(\cub{1})_+ \sm X_2} & & &  \\
{}\save[]*{(\cub{2})_+ \sm X_3} \ar@<1ex>@/^1pc/@{.>}[uurr] \restore & & & &{}\save[]*{(\cub{2})_+ \sm X_2} 
\ar@/_1pc/@{.>}[ul] \restore & \\
& & (1,0,1)}\]
\caption{The object $F_2(f_2,f_3)$.} \label{fil_obj_fig}
\end{figure}
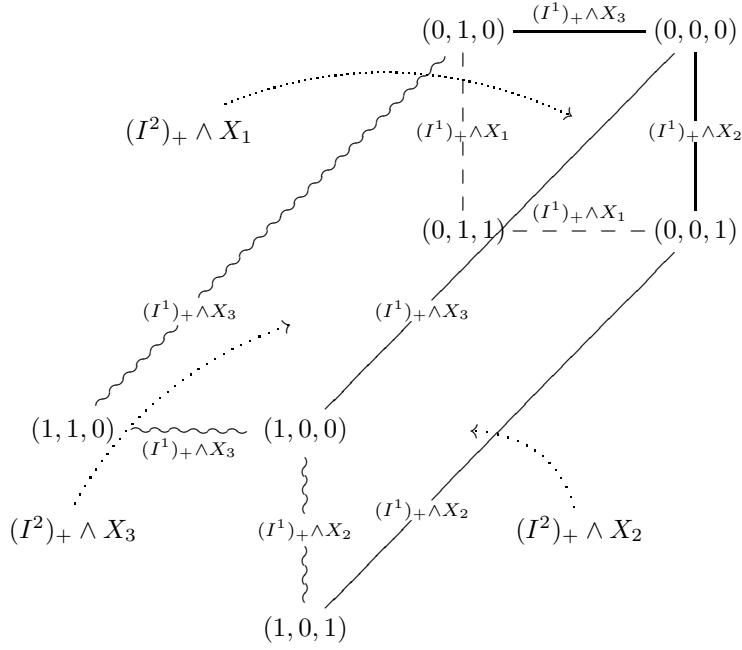
In the diagram, the lines of the shape $\xymatrix@1{~\ar@{~}[r]&~}$ mark the part which is collapsed to the basepoint. 
Thinking of all cubes as subcubes of $\cub{3}$, we glue the 3 objects $(\cub{2})_+ \sm X_1$, $(\cub{2})_+ \sm X_2,$ and
 $(\cub{2})_+ \sm X_3$ (indexed by $(0^i)(\cub{2})$ for $1 \leq i \leq 3$) together along 
two copies of $(\cub{1})_+ \sm X_3$ (indexed by $(0^3)(0^2)(\cub{1})$ and $(0^3)(0^1)(\cub{1})$) and one copy of 
$(\cub{1})_+ \sm X_2$  (indexed by $(0^2)(0^1)(\cub{1})$). Furthermore, we collapse two copies of $(\cub{1})_+ \sm X_3)$ 
(indexed by $(0^3)(1^2)(\cub{1})$ and $(0^3)(1^1)(\cub{1})$) and one copy of $(\cub{1})_+ \sm X_2$ (indexed by $(0^2)(1^1)(\cub{1})$) to the basepoint. 
\end{example}

\begin{lemma} \label{two_maps_giving_tbracket}
The data of the cube system induces maps
\[ \xi_j \colon \randcub{j+1} \sm X_{j+2} \to  F_j(f_2, \dots, f_{j+1}) \quad \text{and} \quad \zeta_j \colon 
F_j(f_2, \dots, f_{j+1}) \to X_0. \]
If $j=n$, the composition $\zeta_n \xi_n$ coincides with $\widetilde{b}^{n+1}(f_1, \dots, f_{n+2})$ of 
Lemma \ref{evalua_descr}. 
\end{lemma}
\begin{proof}
We define $(0^i)^*\xi_j$ to be the composition of  
$(\cub{j})_+ \sm X_{j+2} \to (\cub{j})_+ \sm X_i$ given by the identity smashed with  
$\widetilde{b}^{j+1-i}(f_{i+1}, \dots, f_{j+2})$
using the last $(j+1-i)$ coordinates of the cube and the canonical map $(\cub{j})_+ \sm X_i \to B_j \to F_j(f_2, \dots, f_{j+2})$. 
Its restriction along $1^{k-1}$ is trivial for $k>i$ since $b^{j+1-i}(f_{i+1}, \dots, f_{j+2})$ can be replaced
by the trivial map $b^{j-i}(f_{i+1}, \dots, f_k f_{k+1}, \dots, f_{j+2})$ there. Its restriction along $1^{k}$ is trivial for 
$k < i$ as well, since these subcubes are mapped to the part of $F_j(f_2, \dots, f_{j+1})$ which gets collapsed. 
The same arguments as in Lemma \ref{extending_cube_systems} show that the maps for different $i$ coincide on the
intersections. Therefore, we get an induced
map $\xi_j \colon \randcub{j+1} \sm X_{j+2} \to  F_j(f_2, \dots, f_{j+1})$. 

Next we define $\zeta_j$. On the copy $(\cub{j})_+ \sm X_{i}$ of $B_j$ indexed by $i$ with $1 \leq i \leq j+1$,
we take $\widetilde{b}^{i-1}(f_1, \dots, f_i)$ using the first $(i-1)$ coordinates of $\cub{j}$. This is 
compatible with the identifications of the coequalizer.

To see $\zeta_n \xi_n = \widetilde{c}^{n+1}(f_1, \dots, f_{n+2})$, we look at its restriction to the  subcube
$(0^i)(\cub{n})$. Here $\zeta_n \xi_n$ is  $b^{n+1-i}(f_{i+1}, \dots, f_{n+2})$
using the last $(n+1-i)$ coordinates of the cube, composed with $b^{i-1}(f_1, \dots, f_i)$ using the first $(i-1)$-coordinates.
This is the adjoint of the map which defines $\widehat{b}^{n+1}$ on $(0^i)(\cub{n})$.
\end{proof}

\begin{lemma} 
The object $F_{j+1} := F_{j+1}(f_2, \dots, f_{j+2})$ can be constructed from $F_j:=F_j(f_2, \dots,f_{j+1})$ as 
the mapping cylinder of the map from $F_j$ to the cone $C$ of the map
$\xi_j \colon \randcub{j+1} \sm X_{j+2} \to F_j$. 
The inclusion of $F_j$ into the mapping cylinder therefore gives a map 
$\iota_j \colon F_j \to F_{j+1}$. 
\end{lemma}
\begin{proof}
Let $\tilcub{j+1}$ denote the quotient of $\cub{j+1}$ obtained by collapsing the $(1^i)(\cub{j})$ with $1 \leq i \leq j+1$ to a point. 
Then there is a canonical map 
$\randcub{j+1} \to \tilcub{j+1}$, and we can interpret $\tilcub{j+1}$ as a cone on $\randcub{j+1}$. 
Hence we can model the mapping cone of 
$\xi_j$ by the pushout of
$\tilcub{j+1} \sm X_{j+2} \ot \randcub{j+1} \sm X_{j+2} \xrightarrow{\xi_j} F_{j}.$

To replace the map from $F_j$ to the cone by a cofibration, we need a cylinder object for $F_j$. One
choice for this is $(\cub{1})_+ \sm F_j$, which amounts to adding one additional coordinate to each
$(\cub{i})_+ \sm X_k$ that occurred in the construction of $F_j$. We choose it to be the last coordinate.
Hence the mapping cylinder of $F_j \to C$ 
is weakly equivalent to the pushout of 
\[ \tilcub{j+1} \sm X_{j+2} \ot \randcub{j+1} \sm X_{j+2} \xrightarrow{((0^1)_+ \sm F_j)(\xi_j)} (\cub{1})_+ \sm F_j. \]
The pushout of this diagram is isomorphic to $F_{j+1}$ as defined above. The case $j=1$ can easily be
deduced from Figure \ref{fil_obj_fig}.
\end{proof}

\begin{corollary} \label{Fj_in_triangle}
For $j \leq n$, there is a distinguished triangle
\[ X_{j+2}[j] \xrightarrow{\xi_j} F_j(f_2, \dots, f_{j+1}) \xrightarrow{\iota_{j}} F_{j+1}(f_2, \dots, f_{j+2}) 
\xrightarrow{\pi_{j+1}} X_{j+2}[j+1] \]
in $\Ho(\cC)$. 
\end{corollary}
\begin{proof}
This follows from the last lemma and the definition of the distinguished triangles in the homotopy category of a 
stable model category \cite[Chapter 7]{Hovey-model_1999}.
\end{proof}

\begin{lemma} \label{fil_quotients_fit}
For $0 \leq j \leq n$, the map $ X_{j+2}[j] \xrightarrow{\xi_j}  F_j(f_2, \dots, f_{j+1}) \xrightarrow{\pi_j} X_{j+1}[j]$
equals $f_{j+2}[j]$ in $\Ho(\cC)$.
\end{lemma}
\begin{proof}
The last lemma says that we have a cofibration sequence 
\[F_{j-1}(f_2, \dots, f_j) \xrightarrow{\iota_{j-1}} F_j(f_2, \dots, f_{j+1}) \xrightarrow{\pi_j} X_{j+1}[j]. \]
Hence $\pi_j$ is up to homotopy the map from $F_j$ to its quotient obtained by collapsing every
subcube $(\cub{j})_+ \sm X_i$ of $B_j$ indexed by $2 \leq i \leq j$ to the $(j-1)$-dimensional subcube along which 
it is glued to $(\cub{j})_+ \sm X_{j+1}$. To examine the homotopy class of $\pi_j \xi_j$, we hence only need to
know what $\xi_j$ does on the subcube $(\cub{j})\sm X_{j+2}$ indexed by $j+1$. As it is 
defined to be the map $b^0(f_{j+2})$ on that, we are done.  
\end{proof}

\begin{lemma} \label{F_n_as_fil_obj}
If we consider the $F_j(f_2, \dots, f_{j+1})$ as objects of $\Ho(\cC)$, the sequence  
\[ \ast \to X_1  \xrightarrow{\iota_0} F_1(f_2) 
\xrightarrow{\iota_1} \dots \xrightarrow{\iota_{n-1}} F_n(f_2, \dots, f_{n+1})\]
gives $F_n(f_2, \dots, f_{n+1})$ the structure of an $(n+1)$-filtered object in $\{f_2, \dots, f_{n+1}\}$.  
\end{lemma}
\begin{proof}
We prove that $F_j(f_2, \dots, f_{j+1})$ is a $(j+1)$-filtered object in $\{ f_2, \dots, f_{j+1}\}$ 
by induction.
This is clear for $j=1$. Using that $\pi_j \colon F_j(f_2, \dots ,f_{j+1}) \to X_{j+1}[j]$ plays the role
of the map $\sigma_X$ for $X$ being the $(j+1)$-filtered object 
$F_j(f_2,\dots, f_{j+1})$, we apply Lemma \ref{cones_filtered} 
and Corollary \ref{Fj_in_triangle} to see that $F_{j+1}(f_2, \dots, f_{j+2})$ is a $(j+2)$-filtered object in 
$\{f_2, \dots, f_{j+1}, \pi_j \xi_j[-j]\}$. The last lemma provides the remaining fact $(\pi_j \xi_j) [-j] = f_{j+2}$.
\end{proof}

\begin{proof}[Proof of Proposition \ref{Tbrackets_comp}]
As we have seen in Lemma \ref{two_maps_giving_tbracket}, 
the composition $\zeta_n \xi_n$ is the map $\widetilde{c}^n(f_1, \dots, f_{n+2})$. 
Hence it represents the evaluation of $c$ by Lemma \ref{evalua_descr}.

Let $\sigma_X$ and $\sigma'_X$ denote the structure maps of the filtered object $F_n(f_2, \dots, f_{n+1})$. Lemma
\ref{F_n_as_fil_obj} implies $\sigma_X \xi_n = f_{n}[n-2]$. The definition of $\zeta_n$ and the fact that $\sigma'_X$ 
is the composition 
$ X_1 \xrightarrow{\iota_{(1, 1, \dots, 1)}} (\cub{n+1})_+ \sm X_1 \to F_{n}(f_2, \dots, f_{n+1})$
show $\zeta_n \sigma'_X = f_1$. 
Hence $\xi_n \zeta_n$ is an element of $\tbracket{f_1, \dots, f_{n+2}}$. 
\end{proof}

\section{Relation to $k$-invariants of classifying spaces}

We saw that the evaluation of the universal Toda bracket on a complex is the Toda
bracket of the complex. Since it may as well be evaluated on arbitrary sequences of maps,
it will carry more information than just that about the Toda bracket in general. We will now exhibit how 
its evaluation on a sequence of automorphisms can be expressed. When we
apply our theory to ring spectra in Theorem \ref{gammaR_kinvariants}, this will give us information about the units
of ring spectra (and the units of their matrix rings), rather than only the information about zero divisors 
encoded in the Toda brackets. 

A motivation for this comes from Igusa's results \cite{Igusa_1982} about the first $k$-invariant of the space 
$B\GL_{\infty}(Q \Omega X _+)$, 
which is related to Waldhausen's algebraic $K$-theory of spaces \cite{Waldhausen_1978}: 
Igusa shows that the first $k$-invariant of a connected space $X$ is determined by a cohomology class $k_1^H(\Omega X)$
in the cohomology of the monoid $\pi_0 (\Omega X)$ with coefficients in $H_1(X)$, where the class $k_1^H(\Omega X)$ 
is constructed from the $A_4$-part of the $A_{\infty}$-structure of $\Omega X$ \cite[B, Property 1.1]{Igusa_1982}.
This observation is also used in \cite[Example 4.9, Theorem 3.10]{Baues_D_89}.

We fix a stable topological model category $\cC$, an $n$-split subcategory $\cU$ of $\Ho(\cC)$ 
for some $n \geq 1$, and an $n$-cube system defining $\gamma_{\cU}$. 
We also fix an object $X$ of $\cU$ and denote the representing cofibrant and fibrant object of
$\cC$ which the cube systems chooses as well by $X$.
Consider the topological space $\MapX$ which is pointed by the zero map in $\cC$. Its homotopy groups are 
\[ \pi_i(\MapX,0) \iso [S^i, \MapX]^{\Ho(\Top_*)} \iso [S^i \sm X, X]^{\Ho(\cC)} \iso [X,X]^{\Ho(\cC)}_i. \]
As $\cU$ is $n$-split, $\pi_i(\MapX,0)$ is concentrated
in degrees divisible by $n$. 

The enriched composition in the category $\cC$ equips $\MapX$ with the structure
of a topological monoid, and we refer to the composition as the multiplication. Under the identification above,
the composition of maps in $\Ho(\cC)$ corresponds to the multiplication of $\MapX$. 

The set $\pi_0(\MapX)$ inherits a monoid structure from $\MapX$, and $\MapXinv$ denotes the 
union of all path components of $\MapX$ which are invertible  with respect to the multiplication
on $\pi_0(\MapX)$. Therefore, $\MapXinv$ is a group-like topological monoid. 

As the basepoint of $\MapXinv$ we take $\id_X$, the unit of the multiplication, since the basepoint $0$ of $\MapX$ is not in
$\MapXinv$. There are 
isomorphisms 
\[ \pi_i (\MapX, 0) \iso \pi_i(\MapX, \id_X) \iso \pi_i(\MapXinv, \id_X) \]
for $i \geq 1$.
The second isomorphism is the restriction to the path component.  For the first one, we take the
isomorphism $\sigma_{\id_X}$ of Proposition \ref{discard_basepoint} combined with an adjunction.

A topological monoid $G$ has a classifying space $BG$, defined via the bar construction.
It comes with a map $\omega \colon G \to \Omega BG$.
If the topological monoid $G$ is group-like, that is, the monoid $\pi_0(G)$ is a group, 
then  $\omega$ is a weak equivalence. In our example we get a space $B\MapXinv$ with 
\[ \pi_i(B\MapXinv) \iso \begin{cases} ([X,X]^{\Ho(\cC)})^{\times} & i = 1, \\ 
0 & 1 < i \leq n, \\ [X,X]^{\Ho(\cC)}_n & i=n+1. \end{cases} \]
The left action 
of $\pi_1(B\MapXinv)$ on $\pi_{n+1}(B\MapXinv)$
corresponds under this isomorphism
to the conjugation action $g \cdot \lambda = (g^{-1})^* (g)_* \lambda$ of $[X,X]_0^{\times}$ on $[X[n],X]_0$.

\begin{theorem} \label{toda_br+k-inv}
Let $\cC$ be a stable topological model category, let $\cU$ be a small $n$-split subcategory of $\Ho(\cC)$, 
and let $X$ be a cofibrant and fibrant object of $\cC$ representing an object in $\cU$. Then the restriction map 
\[ \Theta \colon H^{n+2}(\cU, [-,-]_n^{\Ho(\cC)}) \to H^{n+2}(\pi_1(B\MapXinv), \pi_{n+1}(B \MapXinv)) \]
of Proposition \ref{res_map} sends the universal Toda bracket $\gamma_{\cU}$ to the first $k$-invariant
of $B\MapXinv$ not vanishing for dimensional reasons. 
\end{theorem}

We need an auxiliary lemma for the proof. Let $G$ be a group-like topological monoid and let $\omega \colon G \to \Omega BG$
be the map to the group completion. Let $\varphi \colon (S^n,\textrm{pt}) \to (G,g)$ be any map. 
The adjoint of $\omega \varphi$ is
a map from the unreduced suspension $S(S^n)$ to $BG$. It represents an element in $\pi_{n+1}(BG)$. On the other hand, we can 
choose an $h \in G$ such that $gh$ is in the component of $1_G$. If $v$ is a path from $gh$ to $1_G$, we get $[\varphi \cdot h]^v \in \pi_n(G,1_G)$. This does not depend on $v$ and $h
$, as $G$ being a topological monoid implies the $\pi_1$-action on $\pi_n(G)$ to 
be trivial. Composing with $\omega$, we get $\omega_{*}([\varphi \cdot h]^v ) \in   \pi_{n+1} (BG)$.

\begin{lemma} \label{comparing_adjoints}
These two ways to associate an element of $\pi_{n+1}(BG)$ to $\varphi \colon S^n \to G$ are equivalent.
\end{lemma}
\begin{proof}
One can use the homotopy extension property to see that  $S (S^n) \to BG$ is homotopic to a map which sends
$[0,1] \times \{\textrm{pt}\}$ to the basepoint and represents $\omega_{*}([\varphi \cdot h]^v)$.
\end{proof}

\begin{proof}[Proof of Theorem \ref{toda_br+k-inv}] 
We fix a sequence of 
$(n+2)$ automorphisms $(f_1, \dots, f_{n+2})$ of $X$ in $\cU$. Let $f=f_1 \dotsm f_{n+2}$ be their composition. 
We write $\widehat{b}^{n+1}_f$ for  the map $\widehat{b}^{n+1}(f_1, \dots ,f_{n+2}) \colon \skel_n \cub{n+1} \to (\Map_{\cC}(X,X), b^0(f))$.

By the definition of  $c$ in 
Construction \ref{gen_construction} and the restriction
map $\Theta$ in Proposition \ref{res_map}, the evaluation of $\Theta(\gamma_{\cU})$ on $(f_1, \dots, f_{n+2})$ is  
\[ (f^{-1})^* \sigma_{b^0(f)}(\widehat{b}^{n+1}_f) = 
\sigma_{b^0(f)b^0(f^{-1})}((\widehat{b}^{n+1}_f) \cdot b^0(f^{-1})) \in [X,X]_n^{\Ho(\cC)}. \]

We need to examine the image of the homotopy class of this map under  
\[ \begin{split} [X,X]_n^{\Ho(\cC)} & \iso [S^n, (\MapX, 0)] \iso [S^n, (\MapXinv,\id_X)] \\ & \iso [S^{n+1}, B\MapXinv].
\end{split} \]

Choose a path $v$ from $b^0(f) b^0(f^{-1})$ to $\id_X$ in $\Map_{\cC}(X,X)$. Then
\[ \sigma_{\id_X}^{-1} \sigma_{b^0(f)b^0(f^{-1})}[\widehat{b}^{n+1}_f \cdot b^0(f^{-1})] 
=  \sigma_{\id_X}^{-1} \sigma_{\id_X}[\widehat{b}^{n+1}_f \cdot b^0(f^{-1})]^v
= [\widehat{b}^{n+1}_f \cdot b^0(f^{-1})]^v
\]
holds by Proposition \ref{discard_basepoint}, and the last term represents $\Theta(\gamma_{\cU})(f_1, \dots, f_{n+2})$.

For the next step we use Lemma \ref{comparing_adjoints}. It says that the image of 
$[\widehat{b}^{n+1}_f \cdot b^0(f^{-1})]^v$ in $\pi_{n+1}(B\MapXinv)$ is represented by the adjoint of 
$\omega \widehat{b}^{n+1}_f$ considered as a map 
$S (\skel_n \cub{n+1}) \to B\MapXinv$. The unreduced suspension of
$S(\skel_n \cub{n+1})$ is homotopy equivalent to $\randsimpl{n+2}$, the boundary 
of an $(n+2)$-simplex, and we will now explain the resulting map 
$a_f = a(f_1, \dots, f_{n+2}) \colon \randsimpl{n+2} \to B\MapXinv.$

We denote the set of vertices of $\randsimpl{n+2}$ by $\{1, \dots, n+3\}$. Then $a_f$ maps every vertex 
$i$ of $\randsimpl{n+2}$ to the basepoint. The $1$-simplex of $\randsimpl{n+2}$ containing the two
vertices $i < j$ is mapped to $B\MapXinv$ using the path associated to $b^{0}(f_i \dotsm f_{j-1})$ via
the map $\omega \colon \MapXinv \to \Omega B\MapXinv$. Hence every $0$-dimensional subcube of 
$\skel_n \cub{n+1}$ specifies a path from the initial to the terminal vertex of $\randsimpl{n+2}$. This path runs through the 
vertex containing $i<j$ if the term $b^0(f_i \dotsm f_{j-1})$ occurs in the restriction of the cube system
to that $0$-dimensional subcube.

The $2$-simplices of $\randsimpl{n+2}$ containing $i<j<k$ are mapped to $B\MapXinv$ by the 
homotopy between the paths associated to $b^0(f_i \dotsm f_{j-1})$, $b^0(f_j \dotsm f_{k-1})$ and $b^0(f_i \dotsm f_{k-1})$
which we get from $b^1(f_{i} \dotsm f_{j-1}, f_{j} \dotsm f_{k-1})$. This time, the $1$-dimensional subcubes of 
$\skel_n \cub{n+1}$ correspond to the $2$-simplices of $\randsimpl{n+2}$.

The case $n=1$ is displayed in Figure \ref{cube_syst_simplex_fig}, whose right part also appears in 
\cite[B.2.2]{Igusa_1982}.
\begin{figure} 
\begin{minipage}{0.1\textwidth} ~ \end{minipage}%
\begin{minipage}{0.38\textwidth}
\[\xymatrix{ (f)(gh) \ar@{-}[d] \ar@{-}[r] & (f)(g)(h) \ar@{-}[d] \\
(fgh) \ar@{-}[r] & (fg)(h) }\]
\end{minipage}%
\begin{minipage}{0.48\textwidth}
\[\xymatrix{
& & 2 \ar@{-}[dd]^(.7){\omega(g)} \ar@{-}[dddll]_-{\omega(h)} \ar@{-}[dddrr]^-{\omega(fg)} & \\
\\
& & 3 \ar@{-}[dll]_(.25){\omega(gh)} \ar@{-}[drr]^(.25){\omega(f)} \\
1 \ar@{-}[rrrr]_{\omega(fgh)} & & & &  4
}\]
\end{minipage}%
\caption{The square \dots $\qquad \qquad \qquad $  \dots and the associated simplex.} \label{cube_syst_simplex_fig}
\end{figure} 
The situation gets a little bit more involved if $n>1$, since an $(n+1)$-cube has $2(n+1)$ subcubes of dimension 
$n$, but the $(n+2)$-simplex has only $(n+3)$ sub $(n+1)$-simplices. In this case, the $2(n+1)-(n+3)=n-1$ codimension
$1$ subcubes $(0^k)(\cub{n})$ of $\skel_n \cub{n+1}$ with $1 < k < n+1$ do not contribute new information to
the map defined on the boundary of the $(n+2)$-simplex. The reason is that the restriction of the 
pre $(n+1)$-cube 
system to these subcubes is already determined by the underlying $(n-1)$-cube system. We recall that the restriction
to these subcube is built
from $b^{k-1}(f_1, \dots, f_k)$ and $b^{n+1-k}(f_{k+1}, \dots, f_{n+2})$. Accordingly, it corresponds to the restriction 
of the map $a_f \colon \partial \Delta^{n+2} \to B \MapXinv$ to the two simplices with the vertices $\{1, \dots, k\}$
and $\{ k+1, \dots, n+2\}$.  The maps on all other $n$-dimensional subcubes induce maps on one of the $(n+1)$-simplices
of $\partial \Delta^{n+2}$. 

The cochain $(f_1, \dots, f_{n+2}) \mapsto [a(f_1, \dots, f_{n+2})]$ is a representing cocycle for the first $k$-invariant 
as described by Eilenberg and Mac Lane in \cite[\S 19]{Eilenberg_M_1949}. In that reference, the authors also
give an equivalence of this definition of a $k$-invariant to a more commonly used one. 
\end{proof}

\subsection{Coherent vanishing of $k$-invariants} \label{coh_vansh_sec}
The last theorem says  that the vanishing of $\gamma_{\cU}$ implies the vanishing of the first
$k$-invariant of the space $B \MapXinv$ for every cofibrant and fibrant object $X$ of $\cC$ representing an 
object of $\cU$. For our applications, we need a stronger statement in a special case.

For the rest of this section, we assume that $\cC$ is a stable
topological model category in which all objects are fibrant. Furthermore, we assume the $n$-split subcategory
$\cU$ of $\Ho(\cC)$ to have a fixed object $X^1$ such that all other objects
of $\cU$ are finite sums of copies of $X^1$. Such a  $q$-fold sum will
be denoted by $X^q$.

We choose a cofibrant (and automatically fibrant) object of $\cC$ representing $X^1$ and denote it 
also by $X^1$. Let the object $X^q$ in $\cU$  be represented by the $q$-fold coproduct 
$X^1 \wdg \dots \wdg X^1$ of copies of $X^1$ in $\cC$, which we also denote by $X^q$. 
The difference between objects in $\Ho(\cC)$ and $\cC$ will be emphasized by writing $\wdg$ for the 
coproduct in $\cC$ and $\oplus$ for the
coproduct in $\Ho(\cC)$. 

We get maps $\Map_{\cC}(X^q, X^q)  \to \Map_{\cC}(X^{q+1},X^{q+1})$ by adding $\id_{X^1}$ 
on the last summand. 
The restriction of these maps to the set of invertible path components is multiplicative with
respect to the monoid structure. Hence we get a map 
$ t_q \colon B\Map_{\cC}(X^q,X^q)^{\times} \to B\Map_{\cC}(X^{q+1},X^{q+1})^{\times}$ for every $q$.

Here it is convenient to work in a setup with all objects fibrant, since the otherwise necessary fibrant replacement of 
the sum $X^q \wdg X^1$ would mean that we only get a homotopy
class of maps $\Map_{\cC}(X^q, X^q) \to \Map_{\cC}(X^{q+1},X^{q+1})$, rather than an actual map.

Denote the mapping telescope of
$ B\Map_{\cC}(X^1,X^1)^{\times} \to B\Map_{\cC}(X^{2},X^{2})^{\times} \to \dots $
by $B \Map_{\cC}^{\infty}(X,X)^{\times}$. The vanishing of the first $k$-invariant of this space 
does not follow from the vanishing of the first $k$-invariant of all spaces $B\Map_{\cC}(X^q, X^q)^{\times}$
in general, since this  vanishing  does not have to be compatible with the
maps $t_q$. The next lemma provides a sufficient condition for this stronger statement.

\begin{lemma} \label{coh_k-inv_vanish}
Let $\cC$ be a stable topological model category in which all objects are fibrant. Let $\cU$ be a small 
$n$-split subcategory of $\Ho(\cC)$ such that 
\begin{enumerate}[(i)]
\item there is an object $X^1$ in $\cU$ such that all objects of $\cU$ are finite sums of copies of $X^1$,
\item $\gamma_{\cU} \in H^{n+2}(\cU, [-,-]_n^{\Ho(\cC)})$ vanishes,
\item $[X,X]^{\Ho(\cC)}_{i} = 0$ for all objects $X$ of $\,\cU$ if $i>n$, and 
\item $H^{n+1}(\cU, [(-)\oplus X^q, (-) \oplus X^q]_n^{\Ho(\cC)}) = 0$ for all $q \geq 1$. 
\end{enumerate}
Then the space $B \Map_{\cC}^{\infty}(X,X)^{\times}$ has a vanishing $k$-invariant $k^{n+2}$, i.e., it
has the Ei\-len\-berg-Mac Lane space $|B\pi_1(B \Map_{\cC}^{\infty}(X,X)^{\times})|$ as a retract up to 
homotopy. 
\end{lemma}

\begin{proof}
We will construct a section up to homotopy of the $\pi_1$-isomorphism from $B \Map_{\cC}^{\infty}(X,X)^{\times}$
to the Eilenberg-Mac Lane space $|B\pi_1(B \Map_{\cC}^{\infty}(X,X)^{\times})|$. Condition $(iii)$ implies
that it is enough to specify it on $|\skel_{n+2} B\pi_1(B \Map_{\cC}^{\infty}(X,X)^{\times})|$. Since
 $B \Map_{\cC}^{\infty}(X,X)^{\times}$ is constructed as a mapping telescope, it is enough to define 
$\pi_1$-isomorphisms $s_q \colon |\skel_{n+2} B\pi_1(\BMapX{q})| \to \BMapX{q}$ with $t_q s_q \simeq s_{q+1}(t_q)_*$.

As we have seen in the proof of the Theorem \ref{toda_br+k-inv}, the $b^j$ specify maps from all
\mbox{$(j+1)$}-simplices of $|B\pi_1(B \Map_{\cC}^{\infty}(X,X)^{\times})|$ to $\BMapX{q}$. By the compatibility of
the cube system, they assemble to a $\pi_1$-isomorphism $s_q^{n+1}$ defined on the $(n+1)$-skeleton of $|B\pi_1(\BMapX{q})|$

In general, $t_q s_q^{n+1} \simeq s_{q+1}^{n+1}(t_q)_*$ will not hold. But without loss of generality,
we can build this condition into the cube system: for all $j \leq n$, we require  in the
inductive construction of the cube system the map 
$b^j(f_1 \oplus X^1, \dots, f_{j+1} \oplus X^1)$ to be $(- \wdg X^1)b^j(f_1, \dots, f_{j+1})$.

We could extend the $s_q^{n+1}$ to the desired maps $s^q$ if we knew that our cube system extends to
an $(n+1)$-cube system. By Lemma \ref{univ_tbracket_zero}, we know that $\gamma_{\cU} = 0$ implies 
that we can change the maps $b^n(f_1, \dots, f_{n+1})$ to achieve this. Unfortunately, the modified  $b^n$ 
does not have to be compatible with the $t_q$ anymore. 

To fix this, we construct inductively a sequence $e_k \in \ovl{C}^{n+1}(\cU, [-,-]^{\Ho(\cC)}_n)$ for $k \geq 0$ 
with $\delta(e_k) = c$. Since $\gamma_{\cU}=0$, we can find $e_0$ with $\delta(e_0)=c$. Assume we have
built $e_k$ with $\delta(e_k)=c$ and consider 
$e^0_k, e^1_k \in \ovl{C}^{n+1}(\cU, [(-)\oplus X^{k+1} ,(-) \oplus X^{k+1}]^{\Ho(\cC)}_n)$ given by
 \[ e_k^{\epsilon}(f_1, \dots, f_{n+1}) = 
 \begin{cases}
 e(f_1 \oplus X^{k+1}, \dots, f_{n+1} \oplus X^{k+1}) \qquad &\textrm{ if } \epsilon = 0, \\
 (- \wdg X^{k+1})e(f_1, \dots, f_{n+1}) \qquad &\textrm{ if } \epsilon = 1.
 \end{cases}\]
We claim $\delta(e_k^0 - e_k^1) = 0$. This follows from $b^n$ being compatible and 
\[0 = \delta(e_0^{\epsilon})(f_1, \dots, f_{n+2}) + c(f_1 \oplus X^1, \dots, f_{n+1} \oplus X^1) \textrm{ for } \epsilon \in \{1,2\},\]
which holds since the cochain $e_k$ with $\delta(e_k) = c$ can be used to change $b^n$ as in Lemma \ref{univ_tbracket_zero}.
So by (iv), there is an $a_k$ with $\delta(a_k) = e_k^1 - e_k^0$. Now we define $e_{k+1}$ by adding
$\delta(a_k)(f_1, \dots, f_{n+1})$ to $e_k(f_1 \oplus X^{k+1}, \dots, f_{n+1} \oplus X^{k+1})$, and leaving
$e_k$ unchanged on sequences $(f_1, \dots, f_{n+1})$ not of this form. Since $\delta^2(a_k)=0$, we have
$\delta(e_{k+1})=c$. This finishes the construction of the $e_k$. 

We  say that a sequence of 
$(n+1)$ composable maps $(f_1, \dots, f_{n+1})$ in $\cU$ has {\em filtration} $k$ if $k$ is the maximal integer 
such that there exist maps $(f'_1, \dots, f'_{n+1})$
with $f_i = f'_i \oplus X^k$. If we change our compatible $n$-cube system $b^j$ by the cochain $e_{k+1}$ with
the procedure of Lemma \ref{univ_tbracket_zero}, we get an $n$-cube system $b^j_{k+1}$ which is compatible on 
all sequences of filtration up to $k$, and which extends to an $(n+1)$-cube system. Since $b^j_{k+1}$ and $b^j_{k+2}$ 
coincide on the sequences of filtration up to $k$, this is enough to get the desired map on the telescope.
\end{proof}

The hypotheses of the lemma may appear unrealistic at the first glance. Nevertheless, the probably strongest
condition (iv) will reduce to the vanishing of a single Mac Lane cohomology group when we apply it in  
Proposition \ref{splitting_k-theory}. This is much easier to verify as 
to ensure a coherent vanishing of the $k$-invariants by dealing with the associated obstructions on the level of group
cohomology.

\section{The universal Toda bracket of a ring spectrum}      \label{ring_sp_section}
We now apply the results of the preceding sections to ring spectra based on topological spaces. These can be the $S$-algebras of  
\cite{Elmendorf_K_M_M_1997}, the symmetric ring spectra of  \cite{Hovey_S_S_2000} (see \cite{Mandell_M_S_S_2001} for a version
based on topological spaces), or the orthogonal ring spectra introduced in \cite{Mandell_M_S_S_2001}. For a ring spectrum 
$R$, the module category $\modu{R}$ is a topological model category. If $\cC$ is the underlying category of spectra, $ - \sm R \colon \cC \to \modu{R}$ is a left Quillen functor. Hence $\pi_*(R) = [R,R]^{\Ho(\modu{R})}_* \iso [\Sph, R]^{\Ho(\cC)}_*$, and $R$ is
compact in $\Ho(\modu{R})$. 

Recall that $\pi_*(R)$ is $n$-sparse if it is concentrated in degrees divisible by $n$. 
\begin{theorem} \label{univ_higher_tbracket}
Let $R$ be a ring spectrum such that $\pi_*(R)$ is $n$-sparse for some $n \geq 1$. There exists a well defined cohomology class  
$\gamma^{n+2}_R \in \HML^{n+2,-n}_{n\spl}(\pi_*(R))$. By evaluation, it determines
the $(n+2)$-fold Toda bracket of every complex of $(n+2)$ composable maps between finitely generated 
free $n$-sparse $\pi_*(R)$-modules. For a $\pi_*(R)$-module $M$ which admits a resolution by such
modules, the product $\id_M \cup \gamma_R^{n+2}$ is the unique realizability obstruction 
$\kappa_{n+2}(M) \in \Ext^{n+2,-n}_{\pi_*(R)}(M,M)$.  
\end{theorem}

\begin{proof}
Let $\cU$ be the full subcategory of $\Ho(\modu{R})$ given by finite sums of 
copies of the free module of rank $1$ which are shifted by integral multiples of $n$.
Construction \ref{gen_construction} provides 
a cohomology class $\gamma_{\cU} \in H^{n+2}(\cU,[-,-]_{n})^{\Ho(\modu{R})}$. 
The equivalence $\cU \to F(\pi_*(R),n)$ induces an isomorphism 
$H^{n+2}(\cU,[-,-]^{\Ho(\modu{R})}_n) \iso \HML^{n+2,-n}_{n\spl}(\pi_*(R))$. 
The image of $\gamma_{\cU}$ in
the latter group defines $\gamma_{R}$. By Construction \ref{gen_construction} and Theorem \ref{toda_b-obstr},
it has the desired properties.
\end{proof}
We call $\gamma_R^{n+2}$ the universal Toda bracket of $R$. Theorem \ref{triple_universal} is the $n=1$ case of the last theorem. 

\begin{remark} \label{apply_enlargement}
The restriction to modules with a resolution by 
finitely generated free $\pi_*(R)$-modules can be avoided. By \cite[\S 2 and Corollary 3.11]{Jibladze_P_1991}, 
replacing $F(\pi_*(R),n)$ by a larger full small additive 
subcategory doesn't change the cohomology. We choose such a $\cD$ so that it consists only of free modules, and that it contains
all modules from a given free resolution of $M$. Then there is a subcategory $\cU$ in $\Ho(\modu{R})$ equivalent to $\cD$ that gives 
rise to a $\gamma_R$ for which $\id_M \cup \gamma_R$ is defined and equals the obstruction. 
However, there is no small $\cD$ which works for all $M$ simultaneously. Hence
we keep the restriction to the $M$ as stated in the theorem, as this seems to be the most natural choice. \end{remark}

The complex $K$-theory spectrum $KU$ is a ring spectrum \cite[VIII, Theorem 4.2]{Elmendorf_K_M_M_1997} such that 
$\pi_*(KU) \iso \mZ[u^{\pm 1}]$ with $u$ of degree $2$. Its $4$-fold universal Toda bracket is an element of
$\HML^{4,-2}_{2\spl}(\pi_*(KU))$, which is isomorphic to $\HML^4(\mZ)$ by Lemma \ref{HML_gr_ugr} and therefore isomorphic to 
$\mZ/2$ by \cite{Franjou_P_1998}. We compute $\gamma^4_{KU} \neq 0$ in Proposition \ref{computing_gamma_ku}

By \cite[VIII, Theorem 4.2]{Elmendorf_K_M_M_1997} or \cite{Joachim_2001}, the real $K$-theory 
spectrum $KO$ is a ring spectrum. Its graded ring of homotopy groups is given by  
\[\pi_*(KO) = \mZ[ \eta, \omega, \beta^{\pm 1}]/ (2 \eta, \eta^3, \eta \omega, \omega^2 - 4 \beta) \quad \textrm{ with } 
\quad |\eta| = 1, |\omega| = 4, \textrm{ and } |\beta| = 8.\]

The universal Toda bracket $\gamma_{KO} \in \HML^{3,-1}(\pi_*(KO))$ of $KO$ is non-trivial, as $KO$ has non-trivial Toda brackets.
For the reader's convenience, we recall the well known computation  of the easiest example:

\begin{lemma}  
The Toda bracket $\tbracket{2,\eta,2}$ in $\pi_*(KO)$ is defined, has trivial indeterminacy, and contains $\eta^2$.
\end{lemma}
\begin{proof} 
As $2 \eta = 0 = \eta 2$ and $\pi_2(KO)$ is $2$-torsion, the first two statements hold. The ring spectra
map $\Sph \to KO$ is a $\pi_i$-isomorphism for $0 \leq i \leq 2$, so it suffices to calculate the corresponding
Toda bracket for the sphere spectrum. This can be either taken from \cite{Toda-composition_1962} or computed directly, 
following \cite[Theorem 6.1]{Toda_1971}:
Suppose $0 \in \tbracket{2, \eta, 2}$. This would imply the existence of a $4$-cell complex $X$ with $2$, $\eta$ and
$2$ as attaching maps. We consider $H^*(X,\mZ/2)$. Since $\Sq^1$ detects $2$ and $\Sq^2$ detects $\eta$, the existence
of $X$ implies that $\Sq^1 \Sq^2 \Sq^1$ acts non-trivially on the bottom dimensional class in $H^*(X,\mZ/2)$. 
But $\Sq^1 \Sq^2 \Sq^1 = \Sq^2 \Sq^2$, and $\Sq^2\Sq^2$ applied to the bottom class of $H^*(X, \mZ/2)$ is trivial 
for dimensional reasons.
\end{proof}

The class $\gamma_{KO}$ detects non-trivial realizability obstructions:
\begin{lemma} 
The first realizability obstruction $\kappa_3$ of the $\pi_*(KO)$-module $\pi_*(KO) \tensor \mZ/2$ does not vanish. 
Hence $\pi_*(KO) \tensor \mZ/2$ cannot be the homotopy of a $KO$-module spectrum. 

\end{lemma}
\begin{proof}
Write $M$ for $\pi_*(KO) \tensor \mZ/2$. There is a distinguished triangle
$KO \xrightarrow{\cdot 2} KO \to C(2) \to KO[1]$ in $\Ho(\modu{KO})$ which induces a long exact sequence in homotopy. 
Since $M$ is the cokernel of multiplication with $2$ on $\pi_*(KO)$, there is an injection $\iota \colon M \to \pi_*(C(2))$. 
As the two copies of $\pi_*(KO)$ in the long exact sequence are free modules of rank one, $\kappa_3(M)$ vanishes
if and only if $\iota$ split. Hence it is enough to show $\pi_2(C(2)) \iso \mZ/4$. 

From the long exact sequence, we see that $\pi_2(C(2))$ is either $\mZ /4$ or $\mZ /2 \oplus \mZ/2$. 
Let $\rho \in \pi_2(C(2))$ be 
a lift of $\eta \in \pi_1(KO)$ along the epimorphism $\pi_2(C(2)) \to \pi_1(KO)$. Consider

\[ \xymatrix@-0.75pc{ KO[-1] \ar[r] 
& C(2)[-1] \ar[dr] \\
KO[1] \ar@{-->}[u]^{\tau} \ar[r]_{\cdot 2} & KO[1] \ar[u]^{\rho} \ar[r]_{\eta} & KO \ar[r]_{\cdot 2} & KO.} \] 
As we have observed in Remark \ref{different_triple}, a  map $\tau$  such that
the left square commutes is an element of $\tbracket{2,\eta,2}$. 
Hence $2\rho=0$  would imply the contradiction $0 = \tau \in \tbracket{2,\eta,2}$. Therefore, 
$\rho$ cannot be $2$-torsion, and $\pi_2(C(2)) \iso \mZ/4$. 
\end{proof}

\begin{remark} \label{wolberts_thm}
The same argument as in the last lemma shows the corresponding statement about the connective real $K$-theory spectrum. 
This contradicts \cite[Theorem 20]{Wolbert_1998}. 
The reason is an error in \cite[14.1]{Wolbert_1998}. 
In this construction, the author assumes $ku_*$ to be flat as a $ko_*$-module, which does not hold. Accordingly, the generalization
\cite[Theorem 21]{Wolbert_1998} is false as well. 
\end{remark}

\subsection{Universal Toda brackets and $k$-invariants}
Let $R^q$ denote a cofibrant and fibrant object of $\modu{R}$ representing the free $R$-module spectrum of rank $q$. 
We write $\GL_q R$ for the space $\Map_{\modu{R}}(R^q, R^q)^{\times}$
considered in Section \ref{coh_vansh_sec}. 
This definition of the `general linear group' of a ring spectrum $R$ is an
important ingredient for the construction 
of the algebraic $K$-theory of $R$ in the sense of Waldhausen \cite{Waldhausen_1978}, if his definition is interpreted
in the modern language of ring spectra \cite[VI.7]{Elmendorf_K_M_M_1997}. We will encounter the algebraic $K$-theory
of ring spectra in Proposition \ref{splitting_k-theory}.

\begin{theorem} \label{gammaR_kinvariants}
Let $R$ be a ring spectrum such that $\pi_*(R)$ is $n$-sparse for some $n \geq 1$. For $q \geq 1$, the restriction map 
\[ \HML^{n+2,-n}_{n\spl}(\pi_*(R))  \to H^{n+2}(\pi_1 (B \GL_q R), \pi_{n+1} (B \GL_q R))\]
sends the universal Toda bracket $\gamma_R^{n+2}$ of $R$ to the first $k$-invariant of the space $B \GL_q R$
not vanishing for dimensional reasons.
\end{theorem}
\begin{proof}
Since $B\GL_q R = B \Map_{\modu{R}}(R^q, R^q)^{\times}$, this follows from Theorem \ref{toda_br+k-inv} and the
description of the restriction map in Corollary \ref{res_map_HML}.
\end{proof}

Before applying this theorem to examples, we describe the image $\gamma_{R}^{n+2}$ in the ungraded cohomology group
$\HML^{n+2}(\pi_0(R), \pi_n(R))$.

\begin{theorem} \label{ugr_univ_tbracket}
Let $R$ be a ring spectrum with $\pi_*(R)$ concentrated in degrees $0$ and $n$ for some $n \geq 1$. 
There is a universal Toda bracket $\gamma^{n+2}_R \in \HML^{n+2}(\pi_0 (R), \pi_n(R))$. It 
determines all $(n+2)$-fold Toda brackets in $\pi_*(R)$ and the realizability obstruction $\kappa_{n+2}(M)$
of a $\pi_*(R)$-module $M$ which admits a resolution by finitely generated free $n$-sparse $\pi_*(R)$-modules. The 
restriction map 
\[\HML^{n+2}(\pi_0(R), \pi_n(R)) \to H^{n+2}(\pi_1(B\GL_q R), \pi_{n+1}(B \GL_q R))\]
sends $\gamma_R$  to the first $k$-invariant  of 
$B\GL_q R$  not vanishing for dimensional reasons.
\end{theorem}
\begin{proof}
The proof uses the same arguments as that of Theorem \ref{univ_higher_tbracket}. This time,
$\cU$ has the finite sums of (unshifted) copies of $R$ as objects. It is equivalent to $F(\pi_0(R))$. The isomorphism 
$ H^{n+2}(\cU, [-,-]^{\Ho(\cC)}_{n}) \iso \HML^{n+2}(\pi_0(R), \pi_n(R))$ induced by the equivalence enables us to define 
$\gamma_R^{n+2}$ in the latter group. 
\end{proof}

\begin{proposition} \label{res_to_postn_sec}
Let $R$ be a ring spectrum such that $\pi_*(R)$ is $n$-sparse.  
Let $R_{\geq 0}$ be its connective cover and let $P_n(R_{\geq 0})$ be 
the first non-trivial Postnikov section of $R_{\geq 0}$. The restriction
$\HML^{n+2,-n}_{n\spl}(\pi_*(R)) \to \HML^{n+2}(\pi_0(R), \pi_n(R))$
sends the universal Toda bracket of $R$ to the one of $P_n(R_{\geq 0})$. 
\end{proposition}
\begin{proof}
Let $\cU$ be the subcategory of $\Ho(\modu{R})$ given by the finite sums of copies of $R$ which are shifted by 
integral multiples of $n$. The class $\gamma_R^{n+2}$ was defined by applying Construction \ref{gen_construction} to  $\cU$. 
If $\cU_0$ is the subcategory of $\cU$ of finite unshifted copies of $R$, the map from the graded
to the ungraded Mac Lane cohomology is induced by the restriction along the inclusion $\cU_0 \to \cU$. 

Let $\cU_{\geq 0}$ be the subcategory of $\Ho(\modu{R_{\geq 0}})$ which is given by the finite sums 
of unshifted copies of $R_{\geq 0}$. The left Quillen functor
$ - \sm_{R_{\geq 0}} R \colon \modu{R_{\geq 0}} \to \modu{R} $
induces an equivalence between $\cU_{\geq 0}$ and $\cU_0$, since the induced map
on homotopy groups 
\[\modu{\pi_*(R_{\geq 0})} \to \modu{\pi_*(R)}, \quad M \mapsto M \tensor_{\pi_*(R_{\geq 0})} \pi_*(R)\]
restricts to an equivalence between the subcategories of unshifted copies of the free
module of rank $1$. Lemma \ref{morita_cubesys} shows that this equivalence maps the universal Toda bracket of $\cU_0$ 
to the one of $\cU_{\geq 0}$.

A similar argument applied to 
$ - \sm_{R_{\geq 0}} P_n(R_{\geq 0})$
shows that $\gamma_{\cU_0}$ equals the universal Toda bracket of the subcategory of 
$\Ho(P_n R_{\geq 0})$ given by the finite sums of unshifted copies of $P_n R_{\geq 0}$. By Theorem \ref{ugr_univ_tbracket},
this is $\gamma^{n+2}_{P_n R_{\geq 0}}$.
\end{proof} 

We consider the example $KO$ again. The restriction of the universal Toda bracket $\gamma_{KO}$ to $\HML^3(\pi_0(KO), \pi_1(KO))
\iso \HML^{3}(\mZ, \mZ/2)$ is $\gamma_{P_1 KO_{\geq 0}}$. The latter group is $\mZ/2$ \cite[Proposition 13.4.23]{Loday_cyclic_98}.
We show that the image of $\gamma_{KO}$ is the non-zero element, thereby proving once more $\gamma_{KO} \neq 0$. 
Since $P_1 KO_{\geq 0} \iso P_1 ko \iso P_1 \Sph$, this is  a statement about the sphere spectrum, and
computations of Igusa \cite{Igusa_1982} imply
\begin{proposition}
The universal Toda bracket $\gamma_{P_1 \Sph}$ of the first Postnikov section of the sphere spectrum is 
the non-zero element in $\HML^3(\mZ,\mZ/2) \iso \mZ/2$. 
\end{proposition}
\begin{proof}
Let $\cH^m_q$ be the topological monoid of self homotopy equivalences of $q$ copies of the $m$-sphere. Suspension
induces a map $\cH^m_q \to \cH^{m+1}_q$, which is $(m-1)$-connected by the Freudenthal suspension theorem. 

Let $B \cH^m_q$ be the classifying space of $\cH^m_q$. 
The map $\colim_m B\cH^m_q \to B \GL_q \Sph$ is a homotopy equivalence by \cite[Proposition VI.8.3]{Elmendorf_K_M_M_1997}. 
From \cite{Igusa_1982} (compare also \cite[(7.6)]{Baues_D_89}), 
we know that the first $k$-invariant of $B \cH^m_q$
is non-trivial for $q \geq 4$ and $m \geq 3$. The increasing connectivity of the maps in the colimit system
therefore implies that the first $k$-invariant of $B \GL_q \Sph$ does not vanish for $q \geq 4$. Hence
the first $k$-invariant of $B \GL_q P_1 \Sph$ is non-trivial as well. By Theorem \ref{ugr_univ_tbracket}, 
$\gamma_{P_1\Sph}$ has to be non-trivial since the 
$ \HML^3(\mZ,\mZ/2) \to H^{3}(\pi_1 (B \GL_q P_1(\Sph)), \pi_{2} (B \GL_q P_1(\Sph)))$
sends it to this $k$-invariant. 
\end{proof}

Focusing on a ring spectrum with polynomial homotopy again, Proposition \ref{res_to_postn_sec} implies

\begin{corollary} \label{univ_tbracket_gr_vs_ugr}
Let $R$ be a ring spectrum with $\pi_*(R) \iso (\pi_0(R))[u^{\pm 1}]$ for a central unit $u$ in degree $n$. 
The isomorphism $\HML^{n+2,-n}_{n\spl}(\pi_*(R)) \to \HML^{n+2}(\pi_0(R))$ of Lemma \ref{HML_gr_ugr} sends
the universal Toda bracket $\gamma^{n+2}_{R}$ to the one of the first
non-trivial Postnikov section of its connective cover. 
\end{corollary}

This reduces the computation of $\gamma^4_{KU}$ to that of $\gamma^4_{P_2 ku}$.

\begin{remark} \label{rel_to_sp_k-inv}
A ring spectrum $R$ with only two homotopy groups $\pi_0(R)$ and $\pi_n(R)$ has a first $k$-invariant
in the group $\textrm{Der}^{n+1}(\pi_0R, \pi_nR) \iso \THH^{n+2}(\pi_0R,\pi_nR)$ \cite{Lazarev_2001, Dugger_S_2006}. 
Since $\THH^{n+2}(\pi_0R, \pi_nR) \iso \HML^{n+2}(\pi_0R,\pi_nR)$ and the universal Toda brackets coincide 
with the $k$-invariant in the examples $P_1 \Sph$ and $P_2 ku$, we 
expect the universal Toda brackets of first non-trivial Postnikov sections to coincide with these $k$-invariants in general. 
We don't have a proof for this. 
The difficult point is that these two groups are only related by a chain of isomorphisms, 
and we do not know how to identify
the $k$-invariant or the universal Toda bracket in the intermediate steps.

A proof of this statement would not only be interesting for the 
computation of universal Toda brackets. It would also relate the first $k$-invariant 
of a ring spectrum $R$ with the Toda brackets of $R$ and the first $k$-invariants of the spaces $B \GL_q R$ in a very
explicit way. 
\end{remark}

\subsection{A relation to $K$-theory of ring spectra}
For a connective ring spectrum $R$, there is a map $R \to H(\pi_0(R))$ from $R$ to the Eilenberg-Mac Lane
spectrum of $\pi_0(R)$ which is the identity on $\pi_0$. In view of the last remark, we expect the map 
$R \to H(\pi_0(R))$ to split in the homotopy category of ring spectra if $R$ has only two non-trivial homotopy 
groups and a vanishing universal Toda bracket. Though we are not able to prove this statement,
the following proposition will provide a weaker result.

We briefly recall the definition of the algebraic $K$-theory of a ring spectrum $R$, following 
\cite[VI]{Elmendorf_K_M_M_1997}. To avoid technical difficulties, we assume our ring spectrum $R$ to
be an $S$-algebra in the sense of \cite{Elmendorf_K_M_M_1997}. Since all objects in the category of 
$R$-modules are fibrant in this case, 
we obtain maps $B\GL_q R \to B \GL_{q+1}R$ as described in Section \ref{coh_vansh_sec}. 

Let $B\GL R$ be the (homotopy) colimit of the spaces $B\GL_q R$ with respect to these maps.  We apply Quillen's
plus construction to the space $B \GL R$ to obtain $(B \GL R)^+$. For $i \geq 1$, 
algebraic $K$-groups of $R$ can be defined as $K_i(R) = \pi_i((B \GL R)^+)$. We will not need $K_0(R)$, which has
to be defined separately. If $R$ is an Eilenberg-Mac Lane spectrum of a discrete ring $A$, this definition recovers 
the algebraic $K$-groups $K_*(A)$ of $A$ in the sense of Quillen \cite[VI, Theorem 4.3]{Elmendorf_K_M_M_1997}. 

We will later need that the algebraic $K$-theory construction increases connectivity by $1$. 
Recall that map $R \to R'$ of ring spectra is {\em $n$-connected} if the induced map $\pi_i(R) \to \pi_i(R')$
is an isomorphism for $i < n$ and an epimorphism for $i=n$. If $R \to R'$ is $n$-connected, the 
induced map $K_i(R) \to K_i(R')$ is an isomorphism for $i \leq n$ and an epimorphism for $i=n+1$. 
This fact is due to the appearance of the bar construction
in the definition of $K(R)$ and can be proved in a similar way as the corresponding statement
about simplicial rings in \cite[Proposition 1.1]{Waldhausen_1978}. 

\begin{proposition} \label{splitting_k-theory}
Let $R$ be a ring spectrum with homotopy groups concentrated in degrees $0$ and $n$. Suppose that the
universal Toda bracket $\gamma_R^{n+2}$ of $R$ is trivial and that $\HML^{n+1}(\pi_0(R),\pi_n(R))$ 
vanishes. Then the map 
$K_i(R) \to K_i(\pi_0(R))$ induced by $R \to H(\pi_0(R))$ splits for all $i$. 
\end{proposition}
\begin{proof}
It is enough to show that $B\GL R \to B \GL (H(\pi_0(R))$ splits up to homotopy, as this property 
is preserved by the plus construction.  This is equivalent to the splitting of the
map $B \GL R \to | B \pi_1(B \GL R)|$, since both maps are isomorphisms on the fundamental group
and map into an Eilenberg-Mac Lane space. 

We show this by applying Lemma \ref{coh_k-inv_vanish} to the category $\cU$ used in the proof of Theorem \ref{ugr_univ_tbracket}. 
The first three conditions are obviously satisfied. 
It remains to show that $H^{n+1}(\cU, [(-)\oplus R^q, (-)\oplus R^q]_n) = 0$. 

The category $\cU$ is equivalent to $F(\pi_0(R))$. If we
set $A=\pi_0(R)$ and $M=\pi_n(R)$, this equivalence induces an isomorphism
between the last cohomology group and 
\[ H^{n+1}(F(A), \Hom_{A}( (-)\oplus A^q, (- \tensor_{A} M) \oplus M^q)).\]
By Lemma \ref{discard_constant}, this is isomorphic to $\HML^{n+1}(\pi_0(R),\pi_n(R)) \iso 0$.
\end{proof}

\begin{proposition} \label{computing_gamma_ku}
The universal Toda bracket of $KU$ is the non-zero element of $\HML^4(\mZ) \iso \mZ/2$. 
\end{proposition}
\begin{proof}
By Corollary \ref{univ_tbracket_gr_vs_ugr} it is enough to prove  $\gamma^4_{P_2ku} \neq 0$.  
We assume $\gamma^4_{P_2 ku} = 0$ and show that this leads to a contradiction. 

As $\HML^3(\mZ)=0$, Proposition \ref{splitting_k-theory} would imply that $K_3(P_2 ku) \to K_3(\mZ)$ is split. 
This map is onto since $P_2 ku \to H \mZ$ is $2$-connected.
Since $ku \to P_2 ku$ is $4$-connected, $K_3 (ku) \iso K_3(P_2 ku)$, and our assumption implies that $K_3(ku) \to K_3(\mZ)$
is split. 
As the author learned from Ch. Ausoni and J. Rognes, there is a commutative diagram 
\[\xymatrix@-.25pc{
\quad \quad \quad K_3(ku) \ar@{->>}[r] \ar@<3.5ex>[d] & \THH_3(ku) \iso \mZ \quad \ar@{->>}@<-4.5ex>[d] \\
\mZ/48 \iso K _3(\mZ) \ar[r] & \THH_3 (\mZ) \iso \mZ/2
}\]
in which the upper and the right arrow are epimorphisms \cite{Ausoni_2006}. 
Here the horizontal maps are the B{\"o}kstedt trace maps from algebraic $K$-theory to topological Hochschild homology, and 
the vertical maps are induced by $ku \to H(\pi_0(ku))=H(\mZ)$. It follows that the lower map is an epimorphism as well.
If the left map was split, this would mean 
that $\xymatrix@1{\mZ/48 \ar@{->>}[r] & \mZ/2}$ factors through $\mZ$. This is the contradiction implying $\gamma^4_{P_2 ku} \neq 0$.    \end{proof}

One may argue that it would be easier to derive $\gamma^4_{KU} \neq 0$ by calculating  $k^4(B\GL_q P_2 ku) \neq 0$ for some $q \geq 1$.
Unfortunately, we don't know how to compute this $k$-invariant for $q > 1$. For $q=1$, it is trivial. This follows for
example from the restriction $\HML^4(\mZ) \to H^4(\mZ/2, \mZ)$ being trivial. 

\begin{remark}
The $n$th Morava $K$-theory spectrum at a prime $p$ can be represented by a ring spectrum $K(n)$ \cite[\S 11]{Lazarev_2001}.
Here `can be' refers to the fact that there are non-equivalent choices for this structure.

Since $\pi_*(K(n)) \iso \mF_p[v_n^{\pm 1}]$ with $|v_n| = 2(p^n-1)$, we obtain a universal Toda bracket $\gamma^{2p^n}_{K(n)} \in 
\HML^{2p^n}(\mF_p) \iso \mZ/p$. Hence for fixed $p$ and varying $n$, the universal Toda brackets of the $K(n)$ are
elements of $\HML^*(\mF_p)$ lying in the same degrees as the multiplicative generators of the graded ring $\HML^*(\mF_p)$ 
\cite{Franjou_L_S_1994}.

The ring $\pi_*(K(n))$ is a graded field, that is, all $\pi_*(K(n))$-modules are free. Hence it follows that 
all $\pi_*(K(n))$-modules are realizable. Therefore, 
we cannot detect $\gamma^{2p^n}_{K(n)}$ by finding a non-vanishing realizability obstruction.

We do not know if the universal Toda bracket $\gamma^{2p^n}_{K(n)}$ depends on the choice of a model for $K(n)$, and we do also 
not know whether it is non-trivial or not. However, in view of Corollary \ref{univ_tbracket_gr_vs_ugr} and 
Remark \ref{rel_to_sp_k-inv}, we expect $\gamma^{2p^n}_{K(n)}$ to be non-trivial, as the connective Morava $K$-theory
spectrum $k(n)$ has a non-vanishing first $k$-invariant. 
\end{remark}

\begin{remark} 
As mentioned before, Benson, Krause, and Schwede studied a characteristic cohomology class $\gamma_A \in \HH^{3,-1}_k(H^*(A))$
for a differential graded algebra $A$ over a field $k$. Though
they are only concerned with a `triple' characteristic Hochschild class, their theory easily generalizes 
to higher classes when we assume $H^*(A)$ to be $n$-sparse. 
In this case, the Hochschild cochain $m_{n+2}$, which is part of the $A_{\infty}$-structure of $H^*(A)$ \cite{Kadeishvili_1980},
happens to be a Hochschild cocycle. This is easily deduced from the $A_{\infty}$-relations. Similar to the triple class, 
the cohomology class $[m_{n+2}] \in \HH^{n+2,-n}_k(H^*(A))$ is well defined and
determines all $(n+2)$-fold Massey products in $H^*(A)$. 

One may ask whether it is possible to define the higher classes under weaker assumptions, that is, without $H^*(A)$ or
$\pi_*(R)$ being $n$-sparse. To some extent, this is possible if the lower universal classes vanish. We begin by sketching
the first step in the case of a dga $A$. Suppose that $\gamma_A = [m_3] \in \HH^{3,-1}_k(H^*(A))$ vanishes. Then it is possible to 
find an equivalent $A_{\infty}$-structure $(m'_i)$ on $H^*(A)$ such that the {\em cocycle} $m_3$ is zero. This employs  the 
same kind of argument as used to show that every $A_{\infty}$-structure on $H^*(A)$ is trivial if $\HH^{n+2,-n}_k(H^*(A)) = 0$
for $n \geq 1$ \cite{Kadeishvili_1988}. It follows that $m'_4$ is a Hochschild cocycle, which can be used to define a cohomology class
in $\HH^{4,-2}_k(H^*(A))$. This is the candidate for the higher Hochschild class. 
However, it is not unique in general.

In the case of ring spectra, Lemma \ref{univ_tbracket_zero} is the tool for a similar kind of argument. If for example 
$\gamma_R \in \HML^{3,-1}(\pi_*(R))$ vanishes, the lemma says that we can find a $1$-cube system which extends to a $2$-cube system. 
This cube system can be used to define $\gamma_R^4 \in \HML^{4,-2}(\pi_*(R))$ without requiring $\pi_*(R)$ to be $2$-sparse.
As in the algebraic case, there may be different choices for this class $\gamma_R^4$. Moreover, the relation to the Toda brackets
becomes more involved since the indeterminacy is not as easy to control as in the $2$-sparse case. This also affects
the obstruction theory, as there is no unique obstruction class. 
\end{remark}

\appendix

\section{Discarding basepoints of mapping spaces} \label{top_model_cat}
Let $\Top_*$ be the category of pointed compactly generated weak Hausdorff spaces, equipped with the usual model
structure in which the weak equivalences are the weak homotopy equivalences \cite[Theorem 2.4.25]{Hovey-model_1999}. 

A {\em pointed topological model category} $\cC$ is a pointed model category which is enriched, tensored, and cotensored 
over $\Top_*$. This means that there are bifunctors
\[ - \sm -  \colon \Top_* \times \cC \to \cT, \quad 
\Map_{\cC}(-,-)  \colon \cC^{\op} \times \cC \to \Top_*, \quad 
(-)^{(-)}  \colon \cC \times \Top_*^{\op} \to \cC, \]
adjunction isomorphisms 
$
\cC(X, Y^K) \iso \cC(K \sm X, Y) \iso \Top_*(K, \Map_{\cC}(X,Y)), 
$
and an enriched composition $\Map_{\cC}(Y,Z) \sm \Map_{\cC}(X,Y) \to \Map_{\cC}(X,Z).$ The data is asked to satisfy the 
usual associativity and unit conditions. Moreover, the pushout product axiom is required to ensure compatibility of the model 
structures. Details can be found in \cite[4.2]{Hovey-model_1999}.

A {\em stable} topological model category $\cC$ is a pointed topological model category in which 
 the suspension functor $S^1 \sm - \colon \cC \to \cC$
and the loop functor $(-)^{S^1} \colon \cC \to \cC$ form a Quillen equivalence.

For an object $X$ in a category $\cC$, we write $\ucat{X}{\cC}$ for the category of objects under $X$. If $\cC$ is 
a model category, $\ucat{X}{\cC}$ inherits a model structure in which a map is a cofibration, fibration, 
or a weak equivalence if the underlying map in $\cC$ is one \cite[Theorem 7.6.5.(1)]{Hirschhorn-model_2003}.

\begin{proposition} \label{discard_basepoint}
Let $\cC$ be a stable topological model category.
For every
map $g \colon X \to Y$ between cofibrant and fibrant objects in $\cC$, there is an isomorphism 
\[ \sigma_g \colon [S^n, (\Map_{\cC}(X,Y),g)]^{\Ho(\Top_*)} \xrightarrow{\iso} [S^n \sm X, Y]^{\Ho(\cC)}.\]
If $h \colon W \to X$ and $f \colon Y \to Z$ are maps between cofibrant and fibrant objects, the isomorphisms
satisfy $(f_*) (\sigma_g) = (\sigma_{fg})(f_*)$ and $(h^*)(\sigma_g) = (\sigma_{hg})(h^*)$.  For a path $w$ from
$g$ to $g'$ in $\Map_{\cC}(X,Y)$, the isomorphisms $\sigma_g$ and $\sigma_{g'}$ are compatible with the isomorphism 
of homotopy groups induced by $w$, i.e., $\sigma_{g'}(-)^w = \sigma_g$. If $g$ is the zero map, $\sigma_g$ is 
the adjunction isomorphism. 
\end{proposition}

The proof needs some notation and an auxiliary lemma. 
If $S^n$ is the $n$-sphere in $\Top_*$, we consider $S^n_+$ as an object of  $\ucat{S^0}{\Top_*}$. The structure map 
$S^0 \to S^n_+$ sends the basepoint of $S^0$ to the `added' basepoint of $S^n_+$, and the other point to the `original' basepoint
of $S^n$. If $X$ is an object in a pointed topological model category $\cC$, we consider $S^n_+ \sm X$ as an object
in $\ucat{X}{\cC}$ via $X \iso S^0 \sm X \to S^n_+ \sm X$. By the pushout product axiom, $S^n_+ \sm X$ is cofibrant
if $X$ is. If $f \colon X \to Y$ is another object of $\ucat{X}{\cC}$, we write $[S^n_+ \sm X, f]^{\Ho\ucat{X}{\cC}}$
for the set of maps from $X \to S^n_+ \sm X$ to $f$ in the homotopy category of $\ucat{X}{\cC}$.

The following lemma is a reformulation of the well known fact
that $S^n_+$ splits as $S^n \wdg S^0$ after suspension.  
\begin{lemma} \label{discarding_basepoint_in_top}
Let $S^n_+ \sm S^1$ be a space under $S^1$ via $(S^0 \to S^n_+) \sm S^1$, and let $S^{n+1} \wdg S^1$ be a space
under $S^1$ via the inclusion of the second summand. Then there is a map  
$\mu \colon (S^n \sm S^1) \wdg S^1 \to S^n_+ \sm S^1$ in $\ucat{S^1}{\Top}$ which is a homotopy equivalence. 
If $p \colon S^n_+ \to S^n$
is the map which identifies the two basepoints of $S^n_+$,
\[ S^n \sm S^1 \xrightarrow{\textrm{incl}} (S^n \sm S^1) \wdg S^1 \xrightarrow{\mu} S^n_+ \sm S^1 
\xrightarrow{p \sm S^1} S^n \sm S^1 \]
is the identity.
\end{lemma}
\begin{proof}
$S^n$ is a CW-complex with  one $0$-cell and one $n$-cell. The complex 
$S^n_+ \sm S^1 \iso S^n \times S^1 / (S^n \times \{ s_0 \})$ has a $0$-cell, an $1$-cell, and an $(n+1)$-cell. The 
attaching map of the $(n+1)$-cell of $S^n_+ \sm S^1$ is null-homotopic for $n \geq 1$. 
Hence the desired homotopy equivalence exists. If we compose with $p \sm S^1$, we collapse the $1$-cell and do not see
the effect of the null-homotopy. This verifies the last assertion.
\end{proof}

\begin{proof}[Proof of Proposition \ref{discard_basepoint}]
 We define a functor $G \colon \cC \to \cC$ by $G(X) = (X^{S^1})^{\textrm{cof}}$, 
where $(-)^{\textrm{cof}}$ is the functorial cofibrant replacement. The adjunction of suspension and loop gives 
a  natural transformation $\tau \colon S^1 \sm G(X) \to \id_{\cC}$. Since $\cC$ is stable, $\tau_X$ is a weak equivalence
if $X$ is fibrant \cite[Proposition 1.3.13(b)]{Hovey-model_1999}.

Let $\sigma_g$ be the following chain of isomorphisms: 
\[ \begin{split}
&[S^n, (\Map_{\cC}(X,Y),g)]^{\Ho(\Top_*)}  \overset{(i)}{\iso}  
[S^n, (\Map_{\cC}(S^1 \sm G(X),Y),g \tau_X)]^{\Ho(\Top_*)} \\
\overset{(ii)}{\iso} &[S^n_+, (\Map_{\cC}(S^1 \sm G(X),Y),g \tau_X,0)]^{\Ho \ucat{S^0}{\Top_*}} \\ 
 \overset{(iii)}{\iso} &[S^n_+ \sm S^1 \sm G(X), g \tau_X]^{\Ho \ucat{S^1 \sm G(X)}{\cC}} \\
\overset{(iv)}{\iso}&[((S^{n} \sm S^1) \wdg S^1) \sm G(X), g \tau_X]^{\Ho \ucat{S^1 \sm G(X)}{\cC}} \\
\overset{(v)}{\iso}&[S^{n} \sm S^1 \sm G(X), Y]^{\Ho(\cC)}  \overset{(vi)}{\iso} [S^n \sm X, Y]^{\Ho(\cC)}.
\end{split} \] 
Here (i) is induced by the weak equivalence $\tau_X$, (ii) is adding a basepoint, (iii) results from the Quillen 
adjunction between $\ucat{S^0}{\Top_*}$ and $\ucat{X'}{\cC}$ induced by $- \sm X'$ and $\Map_{\cC}(X',-)$ with $X'=S^1 \sm G(X)$, (iv)
uses the weak equivalence under $S^1$ provided by Lemma \ref{discarding_basepoint_in_top}, (v) results from the Quillen 
adjunction between $\ucat{X'}{\cC}$
and $\cC$ given by $Z \mapsto Z \wdg X'$ and the forgetful functor (with $X'=S^1 \sm G(X)$), and (vi) is induced by $\tau_X$ again. 

It is easy to see that the construction is natural in $X$ and $Y$. It is additive since the addition can be 
defined in terms of the $H$-cogroup structure of $S^n$ both in the source and the target. 
Let $w$ be a path from $g$ to $g'$ in $\Map_{\cC}(X,Y)$. Following its action through (i)-(iii), it induces a homotopy
of maps $S^n_+ \sm S^1 \sm G(X) \to Y$ which is itself not a map under $S^1$. But after applying (iv) and (v), the representing
maps become homotopic in $\Ho(\cC)$. 
If $g$ is the zero map, $\sigma_g$ reduces to the adjunction isomorphism: Composing with the $p$ of 
Lemma \ref{discarding_basepoint_in_top} is inverse to (ii), hence Lemma \ref{discarding_basepoint_in_top} shows the assertion.
\end{proof}

\bibliographystyle{abbrv}
\bibliography{univ_tbrs2}
\end{document}